\documentclass[11pt,reqno,a4paper]{amsart}
\usepackage{mathtools,amssymb}
\mathtoolsset{showonlyrefs=true}
\usepackage[foot]{amsaddr}
\usepackage[T1]{fontenc}
\usepackage{lmodern}
\usepackage[expansion=false]{microtype}
\usepackage{mathrsfs,euscript}
\usepackage[only llbracket,rrbracket]{stmaryrd}
\SetSymbolFont{stmry}{bold}{U}{stmry}{m}{n}
\usepackage{tikz}
\usepackage[dvipsnames]{xcolor}
\usepackage[
  final, 
  colorlinks=true,
  linkcolor=MidnightBlue,
  citecolor=MidnightBlue,
  urlcolor=RoyalBlue
]{hyperref}
\hypersetup{
  pdftitle={An elliptic deformation of rationals and continued fractions},
  pdfauthor={Shintarou Yanagida},
  pdfsubject={Elliptic deformation of rationals and continued fractions},
  pdfkeywords={},
  pdfcreator={pdfLaTeX}}

\usepackage{fullpage}
\usepackage{booktabs}
\numberwithin{equation}{section}
\numberwithin{figure}{section}
\allowdisplaybreaks
\usepackage{enumitem}
\setenumerate{label=(\arabic*),nosep}
\setitemize{nosep}

\usepackage{amsthm}
\usepackage[nameinlink,capitalise]{cleveref}
\theoremstyle{definition}
\newtheorem{thm}{Theorem}[section]
\newtheorem{dfn}[thm]{Definition}
\newtheorem{prp}[thm]{Proposition}
\newtheorem{lem}[thm]{Lemma}
\newtheorem{cor}[thm]{Corollary}
\newtheorem{rmk}[thm]{Remark}
\newtheorem{eg}[thm]{Example}
\crefname{thm}{Theorem}{Theorems}
\crefname{dfn}{Definition}{Definitions}
\crefname{prp}{Proposition}{Propositions}
\crefname{lem}{Lemma}{Lemmas}
\crefname{cor}{Corollary}{Corollaries}
\crefname{rmk}{Remark}{Remarks}
\crefname{eg}{Example}{Examples}
\crefname{section}{\S\!}{\S\S\!}
\crefname{subsection}{\S\!}{\S\S\!}
\crefname{equation}{}{}
\crefname{figure}{Figure}{Figures}

\newcommand{\ep}{\epsilon}
\newcommand{\ve}{\varepsilon}
\newcommand{\ol}{\overline}
\newcommand{\wh}{\widehat}
\newcommand{\ot}{\leftarrow}

\newcommand{\lrto}{\leftrightarrow}
\newcommand{\llrto}{\longleftrightarrow}
\newcommand{\mto}{\mapsto}
\newcommand{\lmto}{\longmapsto}
\newcommand{\ceq}{\coloneqq}
\newcommand{\tsum}{{\textstyle \sum}}
\newcommand{\tprd}{{\textstyle \prod}}
\newcommand{\pt}{\mathrm{pt}}
\newcommand{\SL}{\mathrm{SL}}
\newcommand{\bbC}{\mathbb{C}}
\newcommand{\bbQ}{\mathbb{Q}}
\newcommand{\bbZ}{\mathbb{Z}}
\newcommand{\clE}{\mathcal{E}}
\newcommand{\clK}{\mathcal{K}}
\newcommand{\clU}{\mathcal{U}}
\newcommand{\frn}{\mathfrak{n}}
\newcommand{\frG}{\mathfrak{G}}
\newcommand{\scK}{\mathscr{K}}
\newcommand{\sfA}{\mathsf{A}}
\newcommand{\sfB}{\mathsf{B}}
\newcommand{\sfC}{\mathsf{C}}
\newcommand{\sfK}{\mathsf{K}}
\newcommand{\sfM}{\mathsf{M}}
\newcommand{\sfN}{\mathsf{N}}
\newcommand{\sfP}{\mathsf{P}}
\newcommand{\sfT}{\mathsf{T}}
\newcommand{\sfV}{\mathsf{V}}
\newcommand{\bse}{\boldsymbol{e}}
\newcommand{\bsk}{\boldsymbol{k}}
\newcommand{\bsm}{\boldsymbol{m}}
\newcommand{\bsy}{\boldsymbol{y}}
\newcommand{\bsX}{\boldsymbol{X}}
\newcommand{\Cl}{\operatorname{Cl}}
\newcommand{\Gr}{\operatorname{Gr}}
\newcommand{\Ell}{\operatorname{Ell}}
\newcommand{\PSL}{\operatorname{PSL}}
\newcommand{\eell}{\operatorname{e_{\mathrm{ell}}^{\mathrm{mer}}}}
\newcommand{\supp}{\operatorname{supp}}
\newcommand{\dbr}[1]{\llbracket #1 \rrbracket}
\newcommand{\tht}[1]{\theta(#1;p)}
\newcommand{\vth}[1]{\vartheta(#1;p)}
\newcommand{\vc}[2]{\begin{pmatrix}#1\\#2\end{pmatrix}}
\newcommand{\svc}[2]{\begin{psmallmatrix}#1\\#2\end{psmallmatrix}}
\newcommand{\wtf}[2]{w_{a,b;q,p}(#1,#2)}
\newcommand{\rst}[2]{\left. #1 \right|_{#2}}
\newcommand{\bin}[2]{\begin{bmatrix}#1 \\ #2\end{bmatrix}}
\newcommand{\tbin}[2]{\left[\begin{smallmatrix}#1 \\ #2\end{smallmatrix}\right]}
\newcommand{\qbin}[2]{\genfrac{[}{]}{0pt}{}{#1}{#2}_q}

\begin{document}

\title{An elliptic deformation of rational numbers and continued fractions}
\author{Shintarou Yanagida$^1$}
\address{$^1$ Graduate School of Mathematics, Nagoya University,
 Furo-cho, Chikusa-ku, Nagoya, Japan, 464-8602}
\email{yanagida@math.nagoya-u.ac.jp}
\date{November 7, 2026; revised November 27, 2026}

\begin{abstract}
We construct an elliptic deformation $\mathcal{E}_{a,b;q,p}(x)$
of positive rational numbers $x$ by assigning
Schlosser's small elliptic weights $w_{a,b;q,p}(s,t)$
(S\'{e}m.\ Lothar.\ Combin., 2020)
to prefixes of the Stern--Brocot word associated with $x$.
We also introduce elliptic analogues $\mathcal{E}^{\pm}_{a,b;q,p}(x)$
of the regular and negative continued fractions of $x>1$
and show that both coincide with the elliptic rational number
$\mathcal{E}_{a,b;q,p}(x)$.
Under the successive limits $p\to0$, $a\to0$, and $b\to0$,
these elliptic quantities reduce to the $q$-deformed rational numbers
and $q$-continued fractions of Morier-Genoud and Ovsienko
(Forum Math.\ Sigma, 2020).
We further establish elliptic analogues of several formulas and properties
of $q$-rational numbers, including weighted Farey recursions for the numerator
and denominator of $\mathcal{E}_{a,b;q,p}(x)$ and expressions for them
in terms of the closure polynomials (or quiver $F$-polynomials)
of the associated quivers of type $A$.
\end{abstract}

\maketitle
{\small\tableofcontents\par}

\section{Introduction}
\label{s:intro}

Throughout this note, we denote by $\bbZ$ the ring of integers 
and by $\bbQ$ the field of rational numbers.
We also use notation such as $\bbZ_{\ge1}$ and $\bbQ_{>0}$ to denote 
subsets of the integers and rational numbers.

\subsection{Motivation and background}
\label{ss:intro:bg}

Every rational number $x$ satisfying $x>1$ admits 
the regular and negative continued fraction expansions
\begin{align}\label{eq:intro:cf}
  x=a_1+\cfrac{1}{a_2+\cfrac{1}{\ddots+\cfrac{1}{a_{2m}}}}
   =c_1-\cfrac{1}{c_2-\cfrac{1}{\ddots-\cfrac{1}{c_k}}} 
  \qquad 
  (a_i \in \bbZ_{\ge1}, \ c_i \in \bbZ_{\ge 2}), 
\end{align}
denoted by $[a_1,\dotsc,a_{2m}]$ and $\dbr{c_1,\dotsc,c_k}$, respectively.
The classical theory of continued fractions has close connections with 
various mathematical structures, 
including the modular group $\PSL(2,\bbZ)$ \cite[Chaps.\ 1--3]{H},
the Stern--Brocot tree \cite[Chap.\ 14]{Re}, 
and the Farey tessellation \cite{Series}.

Morier-Genoud and Ovsienko \cite{MO} introduced 
a $q$-deformation of positive rational numbers,
and established the compatibility of the $q$-deformed regular and negative 
continued fractions \cite[Theorem~1]{MO}.
For $x=r/s\in\bbQ_{>0}$, their $q$-rational number is denoted by 
\begin{align}\label{eq:intro:xq}
  [x]_q = \frac{r_x(q)}{s_x(q)},
\end{align}
where $r_x(q)$ and $s_x(q)$ have nonnegative integer coefficients.
In the case $x=n$ is a positive integer, we recover the standard $q$-integer
\begin{align}\label{eq:intro:nq}
  [n]_q = 1+q+\dotsb+q^{n-1}.
\end{align}
They also established several $q$-analogues of properties of 
ordinary rational numbers and continued fractions.
In particular, they proved a certain positivity \cite[Theorem~2]{MO}
and weighted Farey recursions \cite[Theorem~3]{MO},
studied the relation to the $q$-deformed modular group \cite[\S4]{MO}.
and interpreted the coefficients of the $q$-rational polynomials 
in terms of the closures of 
the associated type $A$ quivers \cite[Theorem~4]{MO}.

From a different perspective than that of \cite{MO}, 
the relation between continued fractions and cluster combinatorics 
has also been studied using snake graphs in \cite{CS}.
The work of Morier-Genoud and Ovsienko \cite{MO} has 
led to several further developments.
Leclere and Morier-Genoud \cite{LM21} studied $q$-deformations of 
the modular group and real quadratic irrational numbers.
In their subsequent work \cite{LM23}, 
they gave combinatorial interpretations of quantum continuants and rotundi
in terms of paths in triangulated polygons and annuli, respectively.
Another direction was explored by Bapat, Becker, and Licata \cite{BBL}, 
who gave a categorical interpretation of $q$-deformed rational numbers 
using the $2$-Calabi--Yau category of type $A_2$. \par

\medskip

In this note, we investigate an \emph{elliptic deformation} 
of the $q$-rational construction of \cite{MO}.
To explain our motivation, 
let us briefly recall the theory of elliptic hypergeometric functions (HGFs)
(see \cite[Chapter 11]{GR}, \cite{R} and \cite{S}).
In the context of classical and quantum integrable systems, there is 
a well-known hierarchy of rational, trigonometric, and elliptic structures. 
A similar hierarchy appears in the theory of HGFs, consisting of 
ordinary HGFs, basic (or $q$-) HGFs, and elliptic HGFs.
As an example of an elliptic HGF, we recall the very-well-poised 
elliptic hypergeometric series \cite[(11.2.19)]{GR}: 
\begin{align}
\begin{split}
 {}_{r+1}V_{r}&(a_1;a_6,a_7,\dotsc,a_{r+1};q,p;z) \ceq \\
 &\sum_{n=0}^{\infty} \frac{\theta(a_1q^{2n};p)}{\theta(a_1;p)}
  \frac{(a_1,a_6,a_7,\dotsc,a_{r+1};q,p)_n}
       {(q,a_1q/a_6,a_1q/a_7,\dotsc,a_1q/a_{r+1};q,p)_n} (qz)^n.
\end{split}
\end{align}
Without going into the details, we emphasize that 
elliptic HGFs involve theta functions, 
which introduce the \emph{elliptic nome} $p$ in addition to the base $q$.
Motivated by this analogy, we seek a similar elliptic extension of 
the theory of $q$-deformed rational numbers developed in \cite{MO}.

Let us explain the first stage of our approach.
In view of the identity \eqref{eq:intro:nq}, 
it is natural to look for an ``elliptic integer'' 
as a starting point for defining elliptic rational numbers.
To this end, we consider a slightly more general question: 
what would be an elliptic analogue of the Gaussian (or $q$-) 
binomial coefficients $\qbin{n}{k}$?
The theory of elliptic HGFs provides several different approaches 
to this question, including 
Schlosser's elliptic binomial coefficients 
\cite[(2.2)]{Sch0} (see also \cite[\S4]{Sch}) 
and Rosengren's connection coefficients \cite[(2.14)]{R}.

To explain the next stage, let us give some detail of 
the construction of the $q$-rational numbers in \cite{MO}.
Every $x\in\bbQ_{>0}$ is in bijection to the \emph{Stern--Brocot word} 
$U_x=\ve_1\dotsm\ve_N$, where each $\ve_j$ denotes the letter $L$ or $R$ 
(see \cref{ss:w:SB} for the detail).
We consider the matrix product
\begin{align}\label{eq:intro:SB}
  \vc{r_x(q)}{s_x(q)} \ceq B_{\ve_1}(q) \dotsm B_{\ve_N}(q)\vc{1}{1}, \quad
  B_L(q) \ceq \vc{q&0}{q&1}, \quad 
  B_R(q) \ceq \vc{q&1}{0&1}.
\end{align}
Then the entries $r_x(q)$ and $s_x(q)$ give 
the numerator and denominator in \eqref{eq:intro:xq}.
This construction led us to the following idea:
\emph{Replace the base $q$ in the matrices $B_L(q),B_R(q)$ 
 with some elliptic weights},
rather than substitute elliptic integers directly into, say, 
the $q$-continued-fraction formulas.

After exploring several possibilities along these lines, 
we find that \emph{Schlosser's small elliptic weight} \cite[(4.7)]{Sch} 
provides a suitable elliptic analogue of the base $q$.
With the notation introduced in \cref{ss:w:dfn}, 
this weight is given  by
\begin{align}
  w_{a,b;q,p}(s,t) \ceq 
  q\,\frac{\theta(aq^{s+2t},\,bq^{2s+t-2},\,(a/b)q^{t-s-1};p)}
          {\theta(aq^{s+2t-2},\,bq^{2s+t},\,(a/b)q^{t-s+1};p)}
  \quad (s,t\in\bbZ_{\ge1}).
\end{align}
See \eqref{eq:w:sew} for the detail.
In addition to $q$ and $p$, 
the weight depends on two auxiliary parameters $a$ and $b$.
It has the desired limiting properties
\begin{align}\label{eq:intro:wlim}
  \lim_{b\to0}\lim_{a\to0}\lim_{p\to0} w_{a,b;q,p}(s,t) = q, \quad 
  \lim_{q\to1}w_{a,b;q,p}(s,t) = 1.
\end{align}
Hereafter, we call the successive limits $p\to0$, $a\to0$, and $b\to0$ 
the \emph{the $q$-degeneration}.

We describe our main construction and results in \cref{ss:intro:main}.

\subsection{Main constructions and results}
\label{ss:intro:main}

In this subsection, let $a,b,q$ be generic nonzero complex parameters. 
and $p$ be a complex number with $|p|<1$ regarded as the elliptic nome.

Let $w_{a,b;q,p}(s,t)$ denote Schlosser's small elliptic weight \eqref{eq:w:sew}.
For a word $U=\ve_1\dotsm\ve_N$ in the letters $L$ and $R$,
we define its $j$th weight by 
\begin{align}\label{eq:intro:prefix}
  \omega_j(U) \ceq w_{a,b;q,p}\bigl(1+\ell(U_{<j}),1+r(U_{<j})\bigr) \quad 
  (j=1,\dotsc,N), 
\end{align}
where $U_{<j}\ceq\ve_1\dotsm\ve_j$ is the prefix of length $j$,
and $r$ and $\ell$ count its $R$- and $L$-letters (see \eqref{eq:w:omega}).
Thus, an $L$-letter advances the first argument of the weight, 
and an $R$-letter the second.

Our construction of elliptic rational numbers precedes quite similarly to 
the $q$-construction \eqref{eq:intro:SB}, 
using the weights \eqref{eq:intro:prefix} instead of the constant $q$.
Let $x$ be a positive rational number, 
and $U_x=\ve_1\dotsm\ve_N$ be its Stern--Brocot word with $\ve_j\in\{L,R\}$.
Consider the matrix product 
\begin{align}\label{eq:intro:RS}
  \vc{R_x(a,b;q,p)}{S_x(a,b;q,p)} \ceq 
  B_{\ve_1}\bigl(\omega_1(U_x)\bigr) \dotsm B_{\ve_N}\bigl(\omega_N(U_x)\bigr)
  \vc{1}{1},
\end{align}
where the matrices $B_L$ and $B_R$ are the same as in \eqref{eq:intro:SB}.
Now we define the \emph{elliptic rational number} attached to $x$ by 
\begin{align}\label{eq:intro:E}
  \clE_{a,b;q,p}(x) \ceq \frac{R_x(a,b;q,p)}{S_x(a,b;q,p)}.
\end{align}
The first few examples appear in \cref{eg:w:SBcol}.

From the construction \eqref{eq:intro:prefix}--\eqref{eq:intro:E} 
and the property \eqref{eq:intro:wlim}, 
we have  (\cref{prp:w:Elim})
\begin{align}\label{eq:intro:Elim}
  \lim_{b\to0}\lim_{a\to0}\lim_{p\to0} \clE_{a,b;q,p}(x) = [x]_q, \quad 
  \lim_{q\to1}E_{a,b;q,p}(x) = x.
\end{align}
Hence $\clE_{a,b;q,p}(x)$ is a genuine elliptic deformation of 
the $q$-rational $[x]_q$ in \cite{MO}. \par

\medskip

Our first main result identifies this construction with two continued fractions.
Let $x \in \bbQ_{>1}$ in this paragraph.
Using the weights along the extended Stern--Brocot word 
$\wh{U}_x \ceq U_xL$ \eqref{eq:cf:whU},
we define regular and negative \emph{elliptic continued fractions} 
$\clE^{\pm}_{a,b;q,p}(x)$ in \cref{dfn:cf:E+-},
and prove their equality in \cref{thm:cf:clE}:
\begin{align}
  \clE^{\pm}_{a,b;q,p}(x) = \clE_{a,b;q,p}(x).
\end{align}
Under the same limits as in \eqref{eq:intro:wlim} and \eqref{eq:intro:Elim},
the elliptic continued fractions $\clE^{\pm}_{a,b;q,p}(x)$ reduce to 
the $q$-continued fractions of \cite{MO},
and to the ordinary continued fractions in \eqref{eq:intro:cf}.
The further limit $q\to1$ recovers the original rational number $x$.

The second result concerns the \emph{Farey graph},
which will be briefly reviewed in \cref{s:F}.
Recall that two nonnegative reduced rational numbers $x=r/s<y=r'/s'$
are \emph{Farey neighbors} when $r's-rs'=1$,
and that their \emph{mediant} is $x \oplus y \ceq (r+r')/(s+s')$.
Denoting the column vector in \eqref{eq:intro:RS} as 
$V(x) \ceq \svc{R_x(a,b;q,p)}{S_x(a,b;q,p)}$,
We prove the following \emph{weighted Farey recursion} in \cref{thm:F:Farey}: 
\begin{align}
  \sfV(x\oplus y) = \sfV(x)+\frac{\Pi(x\oplus y)}{\Pi(y)}\sfV(y).
\end{align}
Here $\Pi(x)$ is the product \eqref{eq:F:Pi} of the weights along the word of $x$.
In the $q$-degeneration, it reduces to 
the weighted mediant rule underlying \cite[Theorem~3]{MO}.

In addition, \cref{prp:F:inv} proves the elliptic inversion formula
\begin{align}
  \clE_{a,b;q,p}(1/x) = \frac{1}{\clE_{b^{-1},a^{-1};q^{-1},p}(x)}
  \quad (x\in\bbQ_{>0}).
\end{align}
In the $q$-degeneration, it specializes to the inversion formula 
$[1/x]_q=1/[x]_{q^{-1}}$ of the $q$-deformed rational numbers,
which follows from the matrix description of \cite[\S4]{MO}.

There is also a refinement of the two positivity statements in \cite{MO}.
Their Proposition~1.3 establishes positivity for $r_x(q)$ and $s_x(q)$,
while their Theorem~2 establishes coefficientwise positivity of
$r_x(q)s_y(q)-s_x(q)r_y(q)$ for rational numbers $x>y\ge0$.
Now, for our elliptic rational numbers, 
before evaluating the weights as theta quotients,
we may regard all $w_{a,b;q,p}(s,t)$ 
as independent commuting indeterminates $X_{s,t}$.
Both entries of the column then belong to $\bbZ_{\ge0}[\bsX]$ 
with $\bsX\ceq\{X_{s,t}\}_{s,t}$, \cref{prp:F:pos} proves
\begin{align}
  R_x(\bsX)S_y(\bsX)-S_x(\bsX)R_y(\bsX) \in
  \bbZ_{\ge0}[\bsX]\setminus\{0\} \quad (x>y\ge0).
\end{align}

We also establish the elliptic analogue of \cite[Theorem~4]{MO}.
In \cref{s:Q}, to each $x\in\bbQ_{>1}$, 
we attach two type $A$ quivers $Q_x$ and $Q'_x$ 
using the Stern--Brocot word of $x$. 
A subset of vertices of a quiver is called a \emph{closure} 
when every arrow starting in the subset also ends in it.
With the Schlosser weights $y_i$ attached to the vertices 
of $Q_x$ and $Q'_x$, \cref{thm:Q:ellF} gives
\begin{align}\label{eq:intro:F}
  R_x(a,b;q,p) = \sum_{C\in\Cl(Q_x) }\prod_{i\in C}y_i, \quad
  S_x(a,b;q,p) = \sum_{C\in\Cl(Q'_x)}\prod_{i\in C}y_i,
\end{align}
where $\Cl(Q)$ denotes the set of all closures of a quiver $Q$.
The same polynomials are quiver $F$-polynomials for
the thin identity-arrow representations, as discussed in \cref{ss:Q:QG}.

\subsection{Toward a geometric interpretation}
\label{ss:intro:geom}

The main results of this note, outlined in \cref{ss:intro:main}, 
are algebraic in nature.
It is natural to ask whether our elliptic theory admits a geometric 
interpretation, possibly in terms of elliptic cohomology.
Although we do not yet have a satisfactory answer to this question, 
we present a preliminary attempt in this direction in \cref{s:G}.

The vertex weights in \eqref{eq:intro:F} can be written as
normalized elliptic Euler expressions (\cref{s:G}).
Let $E_p \ceq \bbC^\times/p^{\bbZ}$ be an elliptic curve
with elliptic nome $p$ satisfying $0<|p|<1$,
and let $T \ceq (\bbC^\times)^3$ with characters $Y,Z,q$,
so that the $T$-equivariant elliptic cohomology of a point $\pt$
is $\Ell_T(\pt)\simeq E_p^3$.
Consider the virtual $T$-representation \eqref{eq:G:vir}
\begin{align}
  \clK(Y,Z;q) \ceq (1-q^2)(Y+Z-qYZ),
\end{align}
and define its \emph{elliptic Euler expression} by
\begin{align}\label{eq:intro:eell}
  \eell(\clK(Y,Z;q)) \ceq
  q\frac{\theta(Y,Z,q^3YZ;p)}{\theta(q^2Y,q^2Z,qYZ;p)}.
\end{align}
The character balancing makes this quotient 
a meromorphic function on $E_p^3$ (\cref{prp:G:Euler}).
For $Y=bq^{2s+t-2}$ and $Z=(a/b)q^{t-s-1}$,
the expression \eqref{eq:intro:eell} equals $w_{a,b;q,p}(s,t)$.

Along a Stern--Brocot word, an $R$-letter multiplies $(Y,Z)$ by $(q,q)$
and an $L$-letter by $(q^2,q^{-1})$.
In \cref{ss:G:root} we record these shifts by the coweights
$\nu_L\ceq\alpha_1^\vee$ and $\nu_R\ceq\alpha_1^\vee+\alpha_2^\vee$ of $\SL_3$.
Starting with $(Y,Z)=(bq,a/(bq))$, the coweight attached to
the prefix preceding vertex $i$ produces a pair $(Y_i,Z_i)$
with $y_i=\eell(\clK(Y_i,Z_i;q))$ (\cref{prp:G:rwt}).

Combining this with the thin quiver representations of \cref{ss:Q:QG},
whose nonempty quiver Grassmannians are reduced points,
the closure polynomials can be rewritten as finite traces of 
products of these Euler expressions (\cref{prp:G:epush}).

This reformulation follows directly from the definitions.
We record it to make precise which equivariant data an intrinsic
geometric construction would have to produce (see \cref{s:conc}).

\subsection{Organization and conventions}

\Cref{s:w} establishes 
the conventions of theta functions and elliptic weights, 
defines the elliptic Stern--Brocot columns and the elliptic rational numbers,
and proves the integer case formula and the translation formulas.
\Cref{s:cf} introduces regular and negative continued fractions,
and show the compatibility with the elliptic rational numbers.
\Cref{s:F} establishes several properties of elliptic rational numbers,
including weighted Farey recursions, determinant
identities, elliptic inversion and formal positivity.
\Cref{s:Q} proves that the numerator and denominator of
the elliptic rational numbers are equal to the closure polynomials
and quiver $F$-polynomials.
\Cref{s:G} explains a preliminary attempt for geometric interpretation
of our elliptic theory, including the description of the elliptic weights as 
elliptic Euler expressions and the finite-trace reformulation.
\Cref{s:conc} records the remaining limitations and open questions.

Throughout this note, $a,b,q$ denote generic nonzero complex parameters,
and $p$ is a complex number with $0<|p|<1$, called the elliptic nome.
The \emph{$q$-degeneration} means the successive limits
$p\to0$, $a\to0$, and $b\to0$, with $q$ fixed.

\section{Elliptic rational numbers}
\label{s:w}

In this section, we introduce the main object of study in this note,
the \emph{elliptic rational number} $\clE_{a,b;q,p}(x)$ 
associated with each $x\in\bbQ_{>0}$.
We begin by recalling Schlosser's elliptic weights 
and elliptic numbers \cite{Sch}.
For an earlier interpretation in terms of lattice paths, see \cite{Sch0}.

Let $a,b,q$ be generic nonzero complex parameters,
and let $p$ be a complex parameter satisfying $0<|p|<1$.
We regard $p$ as the elliptic nome.

\subsection{Definitions and elementary identities}
\label{ss:w:dfn}

We use the \emph{shifted factorial} \cite[(1.2.19)]{GR}:
\begin{align}\label{eq:w:sf}
  (x;p)_\infty \ceq \prod_{n=0}^\infty (1-xp^n), \quad 
  (x_1,\dotsc,x_d;p)_\infty \ceq \prod_{j=1}^d (x;p)_\infty, 
\end{align}
and the \emph{modified Jacobi theta function} \cite[(11.2.1)]{GR}: 
\begin{align}\label{eq:w:tht}
  \tht{x} \ceq (x;p)_\infty (p/x;p)_\infty, \quad
  \tht{x_1,\dotsc,x_d} \ceq \prod_{j=1}^d \tht{x_j}, 
\end{align}
We also use the \emph{elliptic shifted factorial} %
\cite[(11.2.2)]{GR}, \cite[(2.6)]{W}:
\begin{align}\label{eq:w:esf}
  (x;q,p)_n \ceq \prod_{j=0}^{n-1}\tht{xq^j} \quad (n \in \bbZ_{\ge0}).
\end{align}
We regard the expressions \eqref{eq:w:sf}--\eqref{eq:w:esf} 
as meromorphic functions in $x$ and $x_j$'s.

The fundamental quantity in this note is 
\emph{Schlosser's small elliptic weight} \cite[(4.7)]{Sch}:
\begin{align}\label{eq:w:sew}
  w_{a,b;q,p}(s,t) &\ceq
  q\,\frac{\theta(aq^{s+2t},bq^{2s+t-2},aq^{t-s-1}/b;p)}
          {\theta(aq^{s+2t-2},bq^{2s+t},aq^{t-s+1}/b;p)} \quad 
  (s,t \in \bbZ_{\ge0}).
\end{align}
We regard it as a weight at the lattice point $(s,t)\in\bbZ_{\ge0}^2$.
In the \emph{$q$-degeneration}, i.e., 
the successive limits $p\to0$, $a\to0$, and $b\to0$, we have 
\begin{align}\label{eq:w:wqlim}
  \lim_{b\to0}\lim_{a\to0}\lim_{p\to0} w_{a,b;q,p}(s,t) = q.
\end{align}
On the other hand, in the limit $q\to1$ with $a,b,p$ fixed, we have 
\begin{align}\label{eq:w:wrlim}
  \lim_{q\to1} w_{a,b;q,p}(s,t) = 1.
\end{align}

\begin{eg}\label[eg]{eg:w:weg}
We abbreviate $w_{st}\ceq\wtf{s}{t}$ and record four of these weights explicitly
for later use:
\begin{align}
\begin{split}
  w_{11} &= q\frac{\theta(aq^3,bq,aq^{-1}/b;p)}
                  {\theta(aq,bq^3,aq/b;p)},
  \quad
  w_{12}  = q\frac{\theta(aq^5,bq^2,a/b;p)}
                  {\theta(aq^3,bq^4,aq^2/b;p)},\\
  w_{13} &= q\frac{\theta(aq^7,bq^3,aq/b;p)}
                  {\theta(aq^5,bq^5,aq^3/b;p)},
  \quad
  w_{22}  = q\frac{\theta(aq^6,bq^4,aq^{-1}/b;p)}
                  {\theta(aq^4,bq^6,aq/b;p)}.
\end{split}
\end{align}
\end{eg}

We omit the proof of the following \cref{lem:w:shifts}.

\begin{lem}\label[lem]{lem:w:shifts}
Define two commuting shifts $\sigma_L$ and $\sigma_R$ 
on the parameter pair $(a,b)$ by
\begin{align}
  \sigma_L(a,b)\ceq(aq^2,bq),\quad
  \sigma_R(a,b)\ceq(aq,bq^2).
\end{align}
Then, for $s,t\in\bbZ_{\ge1}$, we have 
\begin{align}\label{eq:w:sigma}
  w_{a,b;q,p}(s+1,t)=w_{\sigma_R(a,b);q,p}(s,t), \quad
  w_{a,b;q,p}(s,t+1)=w_{\sigma_L(a,b);q,p}(s,t).
\end{align}
\end{lem}

Next, let us introduce Schlosser's elliptic numbers.
For that, we denote \emph{Schlosser's elliptic binomial coefficient} %
\cite[(4.4)]{Sch} by 
\begin{align}
  \bin{n}{k}_{a,b;q,p} \ceq
  \frac{(q^{1+k},aq^{1+k},bq^{1+k},aq^{1-k}/b;q,p)_{n-k}}
       {(q,aq,bq^{1+2k},aq/b;q,p)_{n-k}}, \quad
  (k,n\in\bbZ_{\ge0}, \ 0 \le k \le n).
\end{align}
In the $q$-degeneration $p\to0$, $a\to0$, and $b\to0$, we have 
\begin{align}\label{eq:w:binlim}
  \lim_{b\to0}\lim_{a\to0}\lim_{p\to0} \bin{n}{k}_{a,b;q,p} = \bin{n}{k}_q,
\end{align}
where $\qbin{n}{k}$ denotes 
the $q$- (or Gaussian) binomial coefficient \cite[p.~24, 1.2]{GR}.

The special case $n \ge k=1$ is denoted by 
\begin{align}\label{eq:w:Sen}
  [n]_{a,b;q,p} \ceq \bin{n}{1}_{a,b;q,p} =
  \frac{\theta(q^n,aq^n,bq^2,a/b;p)}
       {\theta(q,aq,bq^{n+1},aq^{n-1}/b;p)} \quad (n\in\bbZ_{\ge1}),
\end{align}
and we call it \emph{Schlosser's elliptic number}.
In particular, we have $[1]_{a,b;q,p}=1$ and
\begin{align}
 \lim_{b\to0}\lim_{a\to0}\lim_{p\to0} [n]_{a,b;q,p} = [n]_q,
\end{align}
where $[n]_q \ceq \sum_{j=0}^{n-1}q^j$ is the $q$-number \cite[(1.2.43)]{GR}.

The relationship between the small elliptic weight and the elliptic number 
is given by:

\begin{lem}
For every $n\in\bbZ_{\ge1}$,
\begin{align}\label{eq:w:n-w}
  [n]_{a,b;q,p} = \sum_{t=0}^{n-1} \prod_{j=1}^t w_{a,b;q,p}(1,j).
\end{align}
\end{lem}

\begin{proof}
Denote Schlosser's \emph{big} elliptic weight \cite[(4.6)]{Sch} by
\begin{align}
  W_{a,b;q,p}(s,t) \ceq \prod_{j=1}^t w_{a,b;q,p}(s,j) \quad 
  (s,t\in\bbZ_{\ge0}).
\end{align}
Then Schlosser's elliptic Pascal recursion \cite[(4.5c)]{Sch} gives 
\begin{align}
  \bin{n+1}{k}_{a,b;q,p} = 
  \bin{n}{k}_{a,b;q,p}+W_{a,b;q,p}(k,n+1-k)\bin{n}{k-1}_{a,b;q,p}
\end{align}
for $n,k\in\bbZ$ satisfying $1 \le k \le n$.
Since $\tbin{n}{0}_{a,b;q,p}=1$, in the case $k=1$, it yields 
\begin{align}\label{eq:w:n-W}
  [n+1]_{a,b;q,p}-[n]_{a,b;q,p} = W_{a,b;q,p}(1,n).
\end{align}
Summing \eqref{eq:w:n-W} and using $[1]_{a,b;q,p}=1$, 
we obtain \eqref{eq:w:n-w}.
\end{proof}

\subsection{Elliptic rational numbers}
\label{ss:w:SB}

Finite continued fractions can be encoded by
unimodular $2\times2$ integer matrices.
Morier-Genoud and Ovsienko \cite[\S4.3]{MO} introduced $q$-analogues
of these matrices and studied the corresponding $q$-deformed rational numbers.
Here we construct an elliptic analogue of their argument 
using Schlosser's small elliptic weights \eqref{eq:w:sew}.

Recall the classical Stern--Brocot construction \cite[Chap.~14]{Re}.
We consider the following two transformations of positive rational numbers:
\begin{align}\label{eq:w:LR}
  L\colon x \lmto \frac{x}{x+1}, \quad R\colon x \lmto x+1.
\end{align}
The key observation is that every positive rational number $x$ can be obtained 
from the initial value $1=1/1$ by a unique finite sequence of 
transformations $L$ and $R$.
For example, applying $R$ to $1/1$ gives $2/1$, 
and then applying $L$ gives $2/3$.
Similarly, applying $R$, $L$, and $R$ successively to $1/1$ gives $5/3$.
For each $x\in\bbQ_{>0}$, we record such a sequence of transformations 
by the \emph{Stern--Brocot word} $U_x$ in the letters $L$ and $R$,
written in the order in which the transformations are applied, 
from right to left.
This gives a one-to-one correspondence $x \lrto U_x$  
between positive rational numbers and finite words in the letters $L$ and $R$, 
with the empty word corresponding to $1$.
For example, we have the correspondences
\begin{align}\label{eq:w:SBex}
    R \llrto 2, \quad
   RL \llrto \frac{3}{2}, \quad
  RRL \llrto \frac{5}{2}, \quad
  RLR \llrto \frac{5}{3}.
\end{align}

\begin{figure}[htbp]
\centering
\begin{tikzpicture}[
  x=0.98cm,y=1.1cm,
  every node/.style={inner sep=1.3pt,align=center},
  edgeword/.style={font=\scriptsize,inner sep=0.6pt,fill=white},
  fraction/.style={font=\normalsize} ]
  \node[fraction] (root) at (0,0)
    {$1$\\[-1pt]{\scriptsize $\varnothing$}};
  \node[fraction] (L)  at (-4,-1.25)
    {$\frac{1}{2}$\\[-1pt]{\scriptsize $L$}};
  \node[fraction] (R)  at ( 4,-1.25)
    {$2$\\[-1pt]{\scriptsize $R$}};
  \node[fraction] (LL) at (-6,-2.50)
    {$\frac{1}{3}$\\[-1pt]{\scriptsize $LL$}};
  \node[fraction] (LR) at (-2,-2.50)
    {$\frac{2}{3}$\\[-1pt]{\scriptsize $LR$}};
  \node[fraction] (RL) at ( 2,-2.50)
    {$\frac{3}{2}$\\[-1pt]{\scriptsize $RL$}};
  \node[fraction] (RR) at ( 6,-2.50)
    {$3$\\[-2pt]{\scriptsize $RR$}};
  \node[fraction] (LLL) at (-7,-3.75)
    {$\frac{1}{4}$\\[-1pt]{\scriptsize $LLL$}};
  \node[fraction] (LLR) at (-5,-3.75)
    {$\frac{2}{5}$\\[-1pt]{\scriptsize $LLR$}};
  \node[fraction] (LRL) at (-3,-3.75)
    {$\frac{3}{5}$\\[-1pt]{\scriptsize $LRL$}};
  \node[fraction] (LRR) at (-1,-3.75)
    {$\frac{3}{4}$\\[-1pt]{\scriptsize $LRR$}};
  \node[fraction] (RLL) at ( 1,-3.75)
    {$\frac{4}{3}$\\[-1pt]{\scriptsize $RLL$}};
  \node[fraction] (RLR) at ( 3,-3.75)
    {$\frac{5}{3}$\\[-1pt]{\scriptsize $RLR$}};
  \node[fraction] (RRL) at ( 5,-3.75)
    {$\frac{5}{2}$\\[-1pt]{\scriptsize $RRL$}};
  \node[fraction] (RRR) at ( 7,-3.75)
    {$4$\\[-1pt]{\scriptsize $RRR$}};
  \foreach \parent/\child/\letter in {
    root/L/L,root/R/R,
    L/LL/L,L/LR/R,R/RL/L,R/RR/R,
    LL/LLL/L,LL/LLR/R,LR/LRL/L,LR/LRR/R,
    RL/RLL/L,RL/RLR/R,RR/RRL/L,RR/RRR/R}
    {\draw[thin] (\parent) -- node[edgeword,pos=.54] {$\letter$} (\child);}
\end{tikzpicture}
\caption{\small The Stern--Brocot tree through depth three.
 Each fraction is accompanied by its word, obtained by appending
 $L$ or $R$ along each edge.}
\label{fig:w:SBtree}
\end{figure}
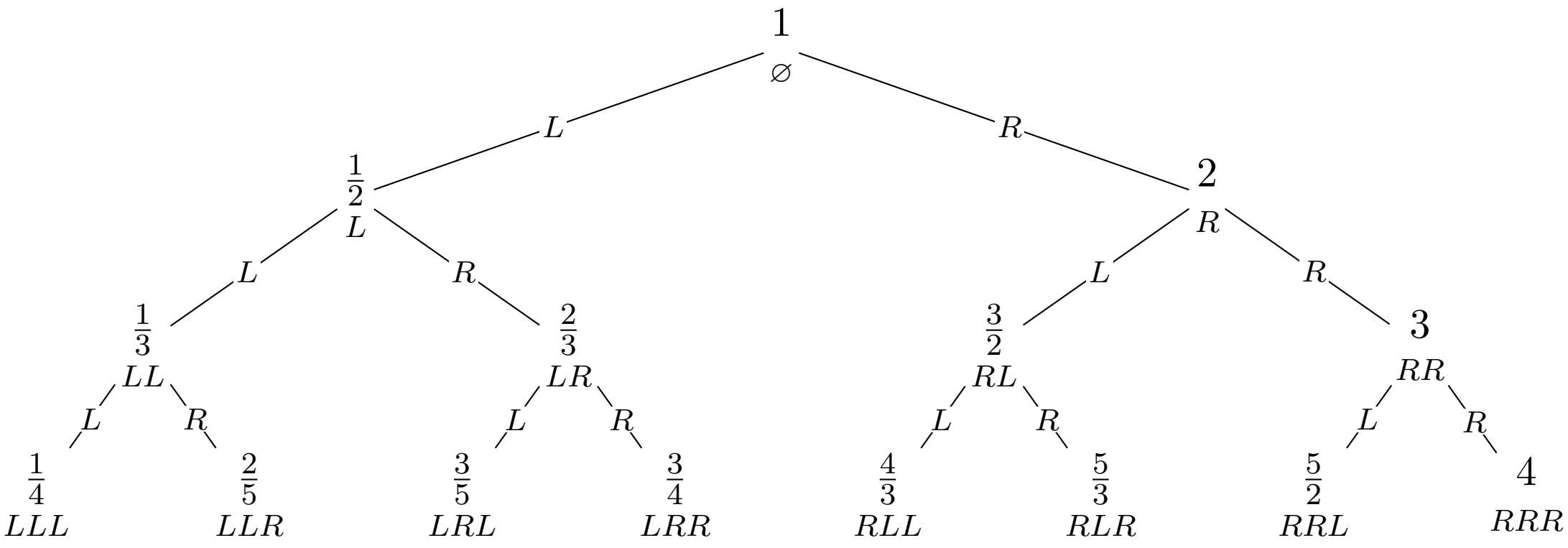

This correspondence can be represented by the \emph{Stern--Brocot tree}.
It is an infinite binary tree whose root corresponds to the empty word, 
or equivalently to the initial value $1/1$.
The two children of a vertex corresponding to a word $W$ are obtained by 
appending $L$ and $R$ to the right of $W$, respectively.
The first three levels of the tree are illustrated in \cref{fig:w:SBtree}.

We can also describe the correspondence $x \lrto U_x$ using the matrix products.
For a scalar $\gamma$, we define $2\times2$ matrices by
\begin{align}\label{eq:w:mat}
  \sfB_L(\gamma) \ceq \begin{pmatrix}\gamma& 0\\\gamma&1\end{pmatrix},
  \quad
  \sfB_R(\gamma) \ceq \begin{pmatrix}\gamma& 1\\0&1\end{pmatrix}.
\end{align}
For a positive rational number $x=r/s$ with $(r,s)$ coprime, 
let $U_x=\ve_1\ve_2\dotsm\ve_N$ be its Stern--Brocot word.
Then, we have 
\begin{align}\label{eq:w:oSB}
  \vc{r}{s} = \sfB_{\ve_1}(1)\sfB_{\ve_2}(1)\dotsm\sfB_{\ve_N}(1)\vc{1}{1},
\end{align}
which follows from the definitions of the transformations \eqref{eq:w:LR}
and the Stern-Brocot word.
See also \cite[\S4.1]{MO}.

\bigskip

Now we will explain the detail of the definition of the elliptic rational number,
which was briefly mentioned in \cref{ss:intro:main}.

For a finite word $U=\ve_1\ve_2\dotsm\ve_N$ in the letters $L$ and $R$,
write $\ell(U)$ and $r(U)$ for the numbers of $L$ and $R$, respectively,
and set $|U| \ceq \ell(U)+r(U) = N$.
For $j=1,2,\dotsc,N$, let $U_{<j}\ceq\ve_1\ve_2\dotsm\ve_{j-1}$ be 
the prefix of length $j$, so that $U_{<1}=\varnothing$.
We put
\begin{align}\label{eq:w:lr}
  \ell_j\ceq\ell(U_{<j}), \quad r_j\ceq r(U_{<j})
\end{align}
whenever the word $U$ is clear from the context. 
Then, using the elliptic weight $w_{a,b;q,p}$ in \eqref{eq:w:sew}, we set
\begin{align}\label{eq:w:omega}
  \omega_j(U) \ceq w_{a,b;q,p}\bigl(1+\ell_j,1+r_j\bigr)
  \quad (j=1,2,\dotsc,N).
\end{align}
Thus $\omega_j(U)$ is the weight at the lattice point
$(1+\ell_j,1+r_j)\in\bbZ_{\ge1}^2$,
reached from $(1,1)$ by increasing the first coordinate by $1$
for each $L$ and the second coordinate by $1$ for each $R$
among the letters preceding $\ve_j$.
Then, for the word $U=\ve_1\ve_2\dotsm\ve_N$, define
\begin{align}\label{eq:w:sfV}
  \sfV(U) \ceq
  \sfB_{\ve_1}\bigl(\omega_1(U)\bigr) \dotsm 
  \sfB_{\ve_N}\bigl(\omega_N(U)\bigr) \vc{1}{1}
\end{align}
with the convention that 
$\sfV(\varnothing) \ceq \svc{1}{1}$ for the empty word $\varnothing$.
Using these notation, we introduce:

\begin{dfn}\label[dfn]{dfn:w:clE}
For each positive rational $x$, let $U_x$ be its Stern--Brocot word.
We define
\begin{align}\label{eq:w:rat}
  \vc{R_x(a,b;q,p)}{S_x(a,b;q,p)} = \sfV_{a,b;q,p}(x) \ceq \sfV(U_x), 
  \quad
  \clE_{a,b;q,p}(x) \ceq \frac{R_x(a,b;q,p)}{S_x(a,b;q,p)}.
\end{align}
We often suppress the dependence on the parameters $a,b,q,p$ and 
simply write $R_x$, $S_x$, $\sfV(x)$, and $\clE(x)$.
We call $\sfV(x)$ and $\clE(x)$ the \emph{elliptic Stern--Brocot column}
and the \emph{elliptic rational number} associated with $x$, respectively.
\end{dfn}

\begin{eg}\label[eg]{eg:w:SBcol}
Using the weights $w_{st}$ of \cref{eg:w:weg},
the words in \eqref{eq:w:SBex} give the following explicit columns.
\begin{align*}
  \sfV(2)&=\sfB_R(w_{11})\vc{1}{1}=\vc{1+w_{11}}{1},\\
  \sfV\left(\frac{3}{2}\right)&=
  \sfB_R(w_{11})\sfB_L(w_{12})\vc{1}{1}=\vc{1+w_{12}+w_{11}w_{12}}{1+w_{12}},\\
  \sfV\left(\frac{5}{2}\right)&=
  \sfB_R(w_{11})\sfB_R(w_{12})\sfB_L(w_{13})\vc{1}{1}
  =\vc{1+w_{11}+w_{13}+w_{11}w_{13}+w_{11}w_{12}w_{13}}{1+w_{13}},\\
  \sfV\left(\frac{5}{3}\right)&=
  \sfB_R(w_{11})\sfB_L(w_{12})\sfB_R(w_{22})\vc{1}{1}
  =\vc{1+w_{12}+w_{12}w_{22}+w_{11}w_{12}+w_{11}w_{12}w_{22}}
      {1+w_{12}+w_{12}w_{22}}.
\end{align*}
In the $q$-degeneration, 
all four weights $w_{st}$ tend to $q$, and these columns specialize to 
the corresponding numerator--denominator pairs of \cite[Example~1.2]{MO}.
\end{eg}

The last observation in \cref{eg:w:SBcol} can be generalized to:

\begin{prp}\label[prp]{prp:w:Elim}
Let $x$ be a positive rational number.
Then, for generic fixed $a,b,p$, we have 
\begin{align}\label{eq:w:Erlim}
  \lim_{q\to1}\clE_{a,b;q,p}(x)=x.
\end{align}
For generic fixed $p$, we have 
\begin{align}\label{eq:w:Eqlim}
  \lim_{b\to0}\lim_{a\to0}\lim_{p\to0}\clE_{a,b;q,p}(x)=[x]_q, 
\end{align}
where $[x]_q$ the $q$-rational of \cite{MO}.
\end{prp}

\begin{proof}
In the limit $q\to1$, every weight $\omega=\omega_j(U_x)$ 
appearing in $\sfV(U)$ tends to $1$ by \eqref{eq:w:wrlim}. 
Thus, \eqref{eq:w:Erlim} follows from \eqref{eq:w:oSB}.

In the $q$-degeneration, every weight tends to $q$ by \eqref{eq:w:wqlim}.
Hence, the matrices $\sfB_R(\omega)$ and $\sfB_L(\omega)$ tend to the 
corresponding matrices in \eqref{eq:intro:SB} (denoted by $R_q$ and $qL_q$,
respectively, in the notation of \cite[Definition~4.4]{MO}).
Now \eqref{eq:w:Eqlim} follows from \eqref{eq:intro:SB} 
(see also \cite[Theorem~1, \S4.6]{MO}).
\end{proof}

As the following \cref{prp:w:trs} shows, 
the orientation of the weight lattice in \eqref{eq:w:omega} is chosen 
so that positive integers give Schlosser's elliptic numbers,
and so that the parameter shifts of \cref{lem:w:shifts}
act along the Stern--Brocot tree.

\begin{prp}\label[prp]{prp:w:trs}
\begin{enumerate}
\item \label{i:prp:w:trs:1}
For every integer $n\ge1$, we have 
\begin{align}\label{eq:w:Eint}
  \sfV(n)=\vc{[n]_{a,b;q,p}}{1}, \quad 
  \clE_{a,b;q,p}(n)=[n]_{a,b;q,p},
\end{align}
where $[n]_{a,b;q,p}$ is Schlosser's elliptic number \eqref{eq:w:Sen}.

\item \label{i:prp:w:trs:2}
For every positive rational $x$,
suppressing the fixed parameters $q,p$ in $\sfV_{a,b;q,p}$, we have
\begin{align}
\label{eq:w:trsV1}
  \sfV_{a,b}\Bigl(\frac{x}{x+1}\Bigr) &= 
  \sfB_L\bigl(w_{a,b;q,p}(1,1)\bigr) \, \sfV_{\sigma_R(a,b)}(x), \\
\label{eq:w:trsV2}
  \sfV_{a,b}(x+1) &= 
  \sfB_R\bigl(w_{a,b;q,p}(1,1)\bigr) \, \sfV_{\sigma_L(a,b)}(x). 
\end{align}
Also, with $w\ceq w_{a,b;q,p}(1,1)$ and 
$\clE_{a,b}(x) \ceq \clE_{a,b;q,p}(x)$, we have
\begin{align}\label{eq:w:trsQ}
  \clE_{a,b}\Bigl(\frac{x}{x+1}\Bigr) =
  \frac{w\,\clE_{\sigma_R(a,b)}(x)}{w\,\clE_{\sigma_R(a,b)}(x)+1}, \quad
  \clE_{a,b}(x+1) = w\,\clE_{\sigma_L(a,b)}(x)+1.
\end{align}
\end{enumerate}
\end{prp}

\begin{proof}
(1)
As can be seen from \eqref{eq:w:LR}, 
the Stern--Brocot word of $n$ is $R^{n-1}$,
and its $j$th letter carries the weight
$\omega_j(R^{n-1})=\wtf{1}{1+j}$.
Since $\sfB_R(\gamma)\svc{X}{1}=\svc{1+\gamma X}{1}$, induction gives
\begin{align}
  \sfB_R(\gamma_1)\dotsm\sfB_R(\gamma_{n-1})\vc{1}{1} =
  \vc{1+\gamma_1+\gamma_1\gamma_2+\dotsb+\gamma_1\dotsm\gamma_{n-1}}{1}.
\end{align}
With $\gamma_j=\wtf{1}{j}$, the top entry is
$\sum_{t=0}^{n-1} \prod_{j=1}^t w_{a,b;q,p}(1,j)=[n]_{a,b;q,p}$
by \eqref{eq:w:n-w}.

(2)
First, note that the maps $x \mto x/(x+1)$ and $x \mto x+1$ are 
the isomorphisms of the Stern--Brocot tree (see \cref{fig:w:SBtree}) 
onto its subtrees rooted at $1/2$ and $2$.
Hence we have $U_{x/(x+1)}=LU_x$ and $U_{x+1}=RU_x$.

To show \eqref{eq:w:trsV1}, 
we consider the weights $\omega_i(U_{x/(x+1)})=\omega_i(LU_x)$ 
for $i=1,2,\dotsc,|U_x|+1$.
We have $\omega_1(LU_x)=\wtf{1}{1}$.
For $1\le j\le|U_x|$, the first identity of \eqref{eq:w:sigma} gives 
\begin{align}
  \omega_{j+1}(LU_x) = \wtf{2+\ell_j}{1+r_j} =
  w_{\sigma_R(a,b);q,p}(1+\ell_j,1+r_j), 
\end{align}
where $r_j \ceq r\bigl((U_x)_{<j}\bigr)$ and 
$\ell_j\ceq\ell\bigl((U_x)_{<j}\bigr)$ (see \eqref{eq:w:lr}).
Hence
\begin{align}
  \sfV_{a,b}\Bigl(\frac{x}{x+1}\Bigr) = 
  \sfB_L(\wtf{1}{1}) \sfV_{\sigma_R(a,b;q,p)}(x),
\end{align}
which is nothing but \eqref{eq:w:trsV1}.

Equation \eqref{eq:w:trsV2} is proved in the same way 
using the second identity of \eqref{eq:w:sigma}.

Finally, \eqref{eq:w:trsQ} follows from
$\sfB_R(w)\svc{X}{Y}=\svc{wX+Y}{Y}$ and
$\sfB_L(w)\svc{X}{Y}=\svc{wX}{wX+Y}$.
\end{proof}

In the $q$-degeneration, \eqref{eq:w:Eint} and \eqref{eq:w:trsQ}
become $[n]_q=1+q+\dotsb+q^{n-1}$ and 
\begin{align}\label{e:qtrs}
  [x+1]_q = q[x]_q+1, \quad 
  \left[\frac{x}{x+1}\right]_q = \frac{q[x]_q}{q[x]_q+1}.
\end{align}
The first identity of \eqref{e:qtrs} is found in \cite[p.~3]{MO}.

\begin{rmk}
We choose the ordered weights in \eqref{eq:w:sfV} so as to obtain
the desired limiting properties (\cref{prp:w:Elim}),
the integer specialization (\cref{prp:w:trs}~\ref{i:prp:w:trs:1}),
and the transformation formulas (\cref{prp:w:trs}~\ref{i:prp:w:trs:2}).
Some of the results presented below depend on 
this particular choice of weights, while others do not.
Let us explain this distinction in more detail.
\begin{enumerate}
\item
The word identities in \cref{s:cf}, the Farey recursions, 
determinant identities, and formal positivity results in \cref{s:F}, 
as well as the quiver closure formulas in \cref{s:Q}, 
are purely algebraic in nature.
They hold for any assignment of weights, provided that 
the weights are assigned consistently to Stern--Brocot prefixes.

\item
In contrast, the inversion formula in \cref{prp:F:inv} relies on 
the reciprocal symmetry \eqref{eq:F:winv} of 
Schlosser's small elliptic weights.
The elliptic Euler expressions in \cref{s:G}
also rely on special properties of these weights.
\end{enumerate}
\end{rmk}

\section{Regular and negative elliptic continued fractions}
\label{s:cf}

In this section, we introduce elliptic analogues of regular and negative
continued fractions, and prove in \cref{thm:cf:clE} that
they are elliptic deformations of the $q$-continued fractions
of Morier-Genoud and Ovsienko \cite[Theorem~1]{MO}.

Recall the matrices $\sfB_R(\gamma)$ and $\sfB_L(\gamma)$ in \eqref{eq:w:mat}.
Additionally, for a nonzero scalar $\gamma$, we define
\begin{align}\label{eq:cf:mat}
  \sfC_L(\gamma) \ceq \gamma^{-1}\sfB_L(\gamma) =
  \begin{pmatrix} 1&0\\1&\gamma^{-1}\end{pmatrix}, \quad
  \sfK \ceq \begin{pmatrix}1&-1\\1&0\end{pmatrix}.
\end{align}
Direct multiplication gives
\begin{align}\label{eq:cf:Kid}
  \sfK\sfB_R(\gamma)=\sfB_L(\gamma)=\gamma\sfC_L(\gamma).
\end{align}

Similarly as in \cref{ss:intro:bg}, following \cite[\S1.1]{MO}, 
we denote the regular and negative continued fraction expressions
of $x\in\bbQ_{>1}$ by
\begin{align}\label{eq:cf:ocf}
 x = [a_1,a_2,\dotsc,a_{2m}] = \dbr{c_1,c_2,\dotsc,c_k}, 
 \quad a_i \ge 1, \ c_i\ge2,
\end{align}
where the regular expansion is assumed to have even length.
In other words,
\begin{align}
 x = a_1+\cfrac{1}{a_2+\cfrac{1}{\ddots+\cfrac{1}{a_{2m}}}}
   = c_1-\cfrac{1}{c_2-\cfrac{1}{\ddots-\cfrac{1}{c_k}}}.
\end{align}
As recalled in \cite[(2.6)]{MO}, the two sequences are related by
\begin{align}\label{eq:cf:c-a}
  (c_1,c_2,\dotsc,c_k) =
  \bigl(a_1+1,2^{a_2-1},a_3+2,2^{a_4-1},\dotsc,a_{2m-1}+2,2^{a_{2m}-1}\bigr),
\end{align}
where $2^j$ denotes a string of $j$ consecutive entries equal to $2$, 
with $2^0$ denoting the empty string.

The Stern--Brocot word $U_x$ of the rational $x$ is then given by
(see \cite[\S4.1]{MO})
\begin{align}\label{eq:cf:URL}
  U_x = R^{a_1} L^{a_2} \dotsm R^{a_{2m-1}} L^{a_{2m}-1},
\end{align}
which has length $|U_x|=\sum_{i=1}^{2m}a_i-1$.
We append one letter $L$ and write
\begin{align}\label{eq:cf:whU}
  \wh{U}_x \ceq U_xL = \ve_1\ve_2\dotsm\ve_M, \quad
  M \ceq |U_x|+1 = \sum_{i=1}^{2m}a_i, \quad
  \ve_M=L.
\end{align}
For $j=1,2,\dotsc,M$, let $\ell_j$ and $r_j$ be 
the letter counts of \cref{ss:w:SB} for the word $\wh{U}_x$.
Now we define
\begin{align}
  v_j \ceq \omega_j(\wh U_x) = \wtf{1+\ell_j}{1+r_j},\quad 1\le j \le M.
\end{align}
In particular, we have $v_j=\omega_j(U_x)$ for $j<M$.

\begin{eg}\label[eg]{eg:vweights}
Let us continue to use the weights $w_{st}$ of \cref{eg:w:weg}
and consider the words in \eqref{eq:w:SBex}.
If $x=5/2$, then $U_x=RRL$ and $\wh U_x=RRLL$, so 
$(v_0,v_1,v_2,v_3) = (w_{11},w_{12},w_{13},w_{23})$.
If $x=5/3$, then $U_x=RLR$ and $\wh U_x=RLRL$, giving
$(v_0,v_1,v_2,v_3) = (w_{11},w_{12},w_{22},w_{23})$.
\end{eg}

For $t\in\bbZ_{\ge0}$ and $n\in\bbZ_{\ge1}$ such that $t+n \le M$, set
\begin{align}
  P_t(n)&\ceq \prod_{j=t}^{t+n-1}v_j = v_t v_{t+1} \dotsm v_{t+n-1},\quad
  P_t(0) \ceq 1,\\
  E_t(n)&\ceq \sum_{i=0}^{n-1}P_t(i)
            = 1+v_t+v_tv_{t+1}+\dotsb+v_t v_{t+1} \dotsm v_{t+n-2},\\
  F_t(n)&\ceq \sum_{i=0}^{n-1}P_t(i)^{-1}
            = 1+\frac{1}{v_t}+\dotsb+\frac{1}{v_t v_{t+1} \dotsm v_{t+n-2}}.
\end{align}
By induction on $n$, we have 
\begin{align}
\label{eq:cf:bl1}
  \sfB_R(v_t)\sfB_R(v_{t+1})\dotsm\sfB_R(v_{t+n-1}) 
  &= \sfM^R_t(n)\ceq\begin{pmatrix}P_t(n)&E_t(n)\\0&1\end{pmatrix},\\
\label{eq:cf:bl2}
  \sfC_L(v_t)\sfC_L(v_{t+1})\dotsm\sfC_L(v_{t+n-1}) 
  &= \sfM^L_t(n)\ceq\begin{pmatrix}1&0\\ F_t(n)&P_t(n)^{-1}\end{pmatrix}.
\end{align}

\begin{dfn}\label[dfn]{dfn:cf:E+-}
Let $x\in\bbQ_{>1}$, and 
consider its regular expansion $[a_1,a_2,\dotsc,a_{2m}]$
and negative expansion $\dbr{c_1,c_2,\dotsc,c_k}$ in \eqref{eq:cf:ocf}.
\begin{enumerate}
\item 
Using $[a_1,a_2,\dotsc,a_{2m}]$, 
we put $t_i \ceq \sum_{h<i}a_h$ for $i=1,2,\dotsc,2m$.
We also set
\begin{align}
  \sfA^+_i \ceq \begin{cases}\sfM^R_{t_i}(a_i) & \text{$i$ odd},\\
                             \sfM^L_{t_i}(a_i) & \text{$i$ even}.\end{cases}
\end{align}
Then we define the \emph{regular elliptic continued fraction} $\clE^+(x)$ by 
\begin{align}
  \clE^+(x) \ceq 
  \text{the first-column ratio of $\sfA^+_1\sfA^+_2\dotsm\sfA^+_{2m}$}.
\end{align}

\item
Using $\dbr{c_1,c_2,\dotsc,c_k}$,
we put $d_i \ceq \sum_{h<i}(c_h-1)$ for $i=1,2,\dotsc,k$.
We also set
\begin{align}\label{eq:cf:A-}
  \sfA^-_i \ceq 
  \sfB_R(v_{d_i})\sfB_R(v_{d_i+1})\dotsm\sfB_R(v_{d_i+c_i-2})\sfK =
  \begin{pmatrix}E_{d_i}(c_i)&-P_{d_i}(c_i-1)\\1&0\end{pmatrix},
\end{align}
using \eqref{eq:cf:bl1} in the last identity.
Then we define the \emph{negative elliptic continued fraction} $\clE^-(x)$ by
\begin{align}
  \clE^-(x) \ceq 
  \text{the first-column ratio of $\sfA^-_1\sfA^-_2\dotsm\sfA^-_k$},
\end{align}
\end{enumerate}
\end{dfn}

\begin{lem}
For $x\in\bbQ_{>1}$ with regular and negative expansions in \eqref{eq:cf:ocf}, 
using the notation in \cref{dfn:cf:E+-}, we have 
\begin{align}\label{eq:cf:E+}
  \clE^+(x) = E_{t_1}(a_1)+
  \cfrac{P_{t_1}(a_1)}{F_{t_2}(a_2)+
  \cfrac{P_{t_2}(a_2)^{-1}}{E_{t_3}(a_3)+
  \cfrac{P_{t_3}(a_3)}{\ddots+
  \cfrac{P_{t_{2m-1}}(a_{2m-1})}{F_{t_{2m}}(a_{2m})}}}}.
\end{align}
and 
\begin{align}\label{eq:cf:E-}
  \clE^-(x) = E_{d_1}(c_1)-
  \cfrac{P_{d_1}(c_1-1)}{E_{d_2}(c_2)-
  \cfrac{P_{d_2}(c_2-1)}{\ddots-
  \cfrac{P_{d_{k-1}}(c_{k-1}-1)}{E_{d_k}(c_k)}}}.
\end{align}
\end{lem}

\begin{proof}
The identity \eqref{eq:cf:E+} can be shown by induction using 
\begin{align}
  \sfM^R_t(n)\vc{X}{Y}=\vc{P_t(n)X+E_t(n)Y}{Y}, \quad
  \sfM^L_t(n)\vc{X}{Y}=\vc{X}{F_t(n)X+P_t(n)^{-1}Y}.
\end{align}
The identity \eqref{eq:cf:E-} follows by induction using
\begin{align}
  \begin{pmatrix}E&-P\\1&0\end{pmatrix}\vc{X}{Y}=\vc{EX-PY}{X}.
\end{align}
\end{proof}

\begin{thm}\label[thm]{thm:cf:clE}
For any $x\in\bbQ_{>1}$, the regular and negative elliptic continued fractions 
coincide with the elliptic rational number defined in \eqref{eq:w:rat}:
\begin{align}
  \clE^+(x) = \clE^-(x) = \clE_{a,b;q,p}(x).
\end{align}
In the $q$-degeneration, we recover the corresponding identities for 
$q$-continued fractions and $q$-rational numbers established by 
Morier-Genoud and Ovsienko \cite[Theorem~1]{MO}.
\end{thm}

\begin{proof}
The classical conversion formula \eqref{eq:cf:c-a} amounts to 
the following identity of words in the letters $R$ and $K$:
\begin{align}\label{eq:cf:wid}
  (R^{c_1-1}K) \dotsm (R^{c_k-1}K)R =
   R^{a_1}(KR)^{a_2} \dotsm R^{a_{2m-1}}(KR)^{a_{2m}}.
\end{align}
Indeed, both sides represent the same word, 
with the letters grouped in different ways.
Setting $M \ceq |U_x|+1=\sum_{i=1}^{2m}a_i$, 
we see that this word contains exactly $M$ occurrences of the letter $R$.

In the word \eqref{eq:cf:wid}, we index the occurrences of the letter $R$
from left to right by $j=1,2,\dotsc,M$.
Then we replace the occurrence indexed by $j$ with 
the matrix $\sfB_R(v_j)$ \eqref{eq:w:mat}
and each occurrence of $K$ with the matrix $\sfK$ \eqref{eq:cf:mat}.
Let $\sfN$ denote the resulting product of matrices.
We compute $\sfN$ in three different ways.

Reading the left-hand grouping of the word \eqref{eq:cf:wid}
and using \eqref{eq:cf:A-}, we obtain
\begin{align}\label{eq:cf:N1}
  \sfN = \sfA^-_1\dotsm\sfA^-_k\sfB_R(v_M).
\end{align}

For the right-hand grouping of \eqref{eq:cf:wid}, 
the description \eqref{eq:cf:URL} shows that
the $j$th letter $R$ lies in a block $R^{a_i}$ with $i$ odd when $\ve_j=R$,
and in a factor $KR$ when $\ve_j=L$.
Now, in the matrix product $\sfN$, 
we replace each $\sfK\sfB_R(v_j)$ by $\sfB_L(v_j)$, 
which does not change $\sfN$ by the first identity in \eqref{eq:cf:Kid}.
Then we get 
\begin{align}\label{eq:cf:N2}
  \sfN = \sfB_{\ve_1}(v_1)\dotsm\sfB_{\ve_M}(v_M).
\end{align}
On the other hand, the second identity in \eqref{eq:cf:Kid} guarantees that
we can replace $\sfK\sfB_R(v_j)$ by $v_j\sfC_L(v_j)$.
Then, using \eqref{eq:cf:bl1} and \eqref{eq:cf:bl2}, we obtain
\begin{align}\label{eq:cf:N3}
  \sfN=\bigl(\tprd_{j: \, \ve_j=L}v_j\bigr)\,\sfA^+_1\dotsm\sfA^+_{2m}.
\end{align}

Now apply $\sfN$ to the column $\svc{1}{0}$.
Since $\sfB_R(\gamma)\svc{1}{0}=\gamma\svc{1}{0}$,
the first description \eqref{eq:cf:N1} gives 
\begin{align}
  \sfN\vc{1}{0} = v_M \cdot 
  (\text{the first column of $\sfA^-_1 \dotsm \sfA^-_k$}).
\end{align}
Since $\ve_M=L$, $v_j=\omega_j(U_x)$ for $j<M$, 
and $\sfB_L(\gamma)\svc{1}{0}=\gamma\svc{1}{0}$,
the second \eqref{eq:cf:N2} gives 
\begin{align}
  \sfN\vc{1}{0} = v_M \cdot \sfV(U_x).
\end{align}
Finally, the third description \eqref{eq:cf:N3} gives 
\begin{align}
  \sfN\vc{1}{0} = (\text{a non-zero scalar}) \cdot 
  (\text{the first column of $\sfA^+_1 \dotsm \sfA^+_{2m}$}).
\end{align}
Hence the three columns are proportional,
and we have $\clE^-(x)=\clE_{a,b;q,p}(x)=\clE^+(x)$.

In the $q$-degeneration, all the weights $v_j$ tend to $q$,
so \eqref{eq:cf:E-} and \eqref{eq:cf:E+} become 
the $q$-continued fractions of \cite[Theorem~1]{MO}, 
and their common value is $[x]_q$ by \cref{prp:w:Elim}.
\end{proof}

\begin{eg}
We check \cref{thm:cf:clE} for the four cases in \cref{eg:w:SBcol}.
Let us write $\clE(x) \ceq \clE_{a,b;q,p}(x)$ 
and use the weights $w_{st}$ of \cref{eg:w:weg}.

For $x=2=[1,1]=\dbr{2}$, we have $\clE(2)=1+w_{11}$ and
\[
  \clE^+(2) = E_0(1)+\frac{P_0(1)}{F_1(1)} = 1+w_{11},\quad
  \clE^-(2) = E_0(2)=1+w_{11}.
\]

For $x=3/2=[1,2]=\dbr{2,2}$,
we have $\clE(3/2)=(1+w_{12}+w_{11}w_{12})/(1+w_{12})$ and
\[
  \clE^+\!\left(\frac{3}{2}\right) 
  =E_0(1)+\frac{P_0(1)}{F_1(2)}
  =1+\frac{w_{11}}{1+w_{12}^{-1}},\quad
  \clE^-\!\left(\frac{3}{2}\right)
  =E_0(2)-\frac{P_0(1)}{E_1(2)}
  =1+w_{11}-\frac{w_{11}}{1+w_{12}}.
\]

For $x=5/2=[2,2]=\dbr{3,2}$, we have
\begin{align*}
  \clE\!\left(\frac{5}{2}\right)
  &=\frac{1+w_{11}+w_{13}+w_{11}w_{13}+w_{11}w_{12}w_{13}}{1+w_{13}}, \\
  \clE^+\!\left(\frac{5}{2}\right)
  &=E_0(2)+\frac{P_0(2)}{F_2(2)}
   =1+w_{11}+\frac{w_{11}w_{12}}{1+w_{13}^{-1}},\\
  \clE^-\!\left(\frac{5}{2}\right)
  &=E_0(3)-\frac{P_0(2)}{E_2(2)}
   =1+w_{11}+w_{11}w_{12}-\frac{w_{11}w_{12}}{1+w_{13}}.
\end{align*}

Finally, for $x=5/3=[1,1,1,1]=\dbr{2,3}$, we obtain
\begin{align*}
  \clE\!\left(\frac{5}{3}\right)
  &=\frac{1+w_{12}+w_{12}w_{22}+w_{11}w_{12}+w_{11}w_{12}w_{22}}
         {1+w_{12}+w_{12}w_{22}}, \\
  \clE^+\!\left(\frac{5}{3}\right)
  &=E_0(1)+\cfrac{P_0(1)}{F_1(1)+
           \cfrac{P_1(1)^{-1}}{E_2(1)+\cfrac{P_2(1)}{F_3(1)}}}
   =1+\cfrac{w_{11}}{1+\cfrac{w_{12}^{-1}}{1+w_{22}}},\\
  \clE^-\!\left(\frac{5}{3}\right)
  &=E_0(2)-\frac{P_0(1)}{E_1(3)}
   =1+w_{11}-\frac{w_{11}}{1+w_{12}+w_{12}w_{22}}.
\end{align*}
In each case, the two nested elliptic fractions $\clE^{\pm}(x)$
and the ratio $\clE(x)$ of the columns in \cref{eg:w:SBcol} agree.
The weight $v_M$ of the appended letter does not appear in these ratios.
\end{eg}

\section{Farey recursions, inversion, and determinant identities}
\label{s:F}

In this section, we study how the elliptic Stern--Brocot columns %
$\sfV(x)$ (\cref{dfn:w:clE}) behave under Farey mediation, 
inversion and determinant formation.
The formulas obtained in this section refine the corresponding
results of \cite[Theorems 2 and 3]{MO}.

Let us begin by recalling the classical \emph{Farey graph}.
For further details, we refer to \cite[Chap.~1]{H}, 
where it is called the Farey diagram.
Its vertex set is $\bbQ\cup\{\infty\}$, where $\infty=1/0$.
Two reduced fractions $r/s$ and $r'/s'$ with nonnegative denominators
are joined by an edge if $|r's-rs'|=1$.
Two vertices joined by an edge are called \emph{Farey neighbors}.
In particular, $0=0/1$ and $\infty=1/0$ are Farey neighbors, 
and we call the edge $(0,\infty)$ the \emph{initial edge}.

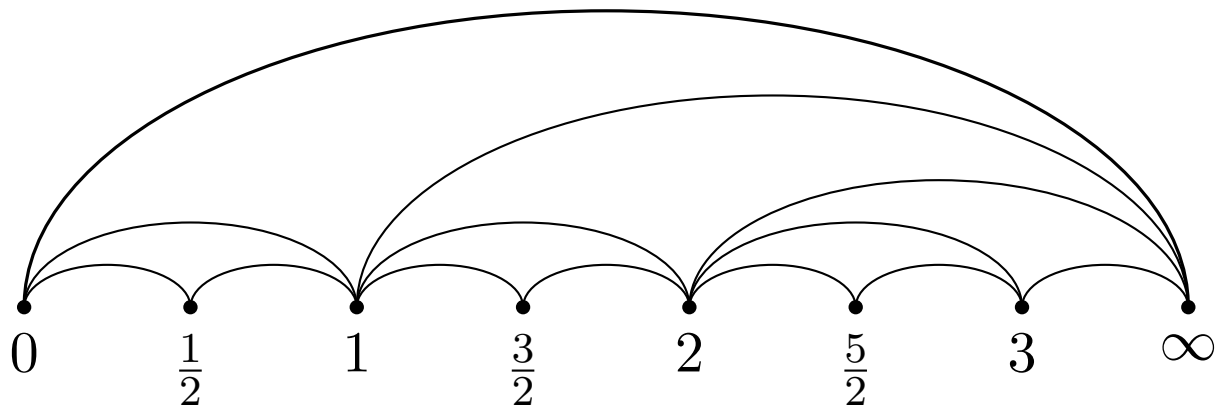
\begin{figure}[htbp]
\centering
\begin{tikzpicture}[x=1cm,y=1cm,>=stealth,
  every node/.style={inner sep=1pt}]
  \foreach \i/\j in {0/1,1/2,0/2,2/3,3/4,2/4,4/5,5/6,4/6,
                        2/7,4/7,6/7} {
    \draw[thin] ({1.1*\i},0)
       arc[start angle=180,end angle=0,
           x radius={0.55*(\j-\i)},y radius={0.28*(\j-\i)}];
  }
  \draw[semithick] (0,0)
     arc[start angle=180,end angle=0,x radius=3.85,y radius=1.96];
  \foreach \i/\t in {0/0,1/{\frac{1}{2}},2/1,3/{\frac{3}{2}},
                        4/2,5/{\frac{5}{2}},6/3,7/\infty} {
    \fill ({1.1*\i},0) circle (1.35pt);
    \node[below=4pt] at ({1.1*\i},0) {$\t$};
  }
\end{tikzpicture}
\caption{\small A finite subgraph of the Farey graph. 
 The outer arc is the initial edge $(0,\infty)$, and $\frac{3}{2}=1\oplus2$.}
\label{fig:F:Farey}
\end{figure}

Starting with the initial edge, we first insert the \emph{mediant}
$(0/1)\oplus(1/0)\ceq(0+1)/(1+0)=1/1$
and then repeatedly insert the mediant $(r/s)\oplus(r'/s')\ceq(r+r')/(s+s')$
of each pair $(r/s,r'/s')$ of consecutive vertices.
Every positive rational number appears exactly once in this process, 
and the two vertices between which it is inserted 
are called its \emph{Farey parents}.
It is a classical fact that 
every pair of Farey neighbors $0\le x<y\le\infty$ appears 
as a pair of consecutive vertices at some stage of this process.
Thus, $x$ and $y$ are the Farey parents of their mediant $x\oplus y$.
A finite subgraph of the Farey graph is illustrated in \cref{fig:F:Farey}.

To formulate the recursion relations for the elliptic Stern--Brocot columns 
$\sfV(x)$ in terms of the Farey graph, 
we define the \emph{accumulated weight} of a word $U$ 
in the letters $L$ and $R$ by
\begin{align}\label{eq:F:Pi}
  \Pi(U) \ceq \prod_{j=1}^{|U|}\omega_j(U), \quad
  \Pi(\varnothing) \ceq 1,
\end{align}
where $\omega_j(U)$ is defined in \eqref{eq:w:omega}.
For each $x\in\bbQ_{>0}$, let $U_x$ be its Stern--Brocot word and 
set $\Pi(x)\ceq\Pi(U_x)$.
We also set $\Pi(\infty)\ceq1$ by convention.

\begin{thm}\label{thm:F:Farey}
Let $0\le x=r/s<y=r'/s'\le\infty$ be Farey neighbors satisfying
$r's-rs'=1$, and let $z\ceq x\oplus y=(r+r')/(s+s')$ be their mediant.
Then we have
\begin{align}\label{eq:F:Farey}
  \sfV(z)=\sfV(x)+\frac{\Pi(z)}{\Pi(y)}\sfV(y).
\end{align}
\end{thm}

\begin{proof}
We examine the Stern--Brocot words $U_z,U_x,U_y$ of 
the mediant $z$ and of its Farey parents $x<y$.
All the possible patterns of these words can be divided into 
the five cases of the table below.
The first two rows describe the case 
when $U_z$ contains both letters $L$ and $R$, and 
the next two rows describe the case 
when $U_z$ consists of only one kind of letter.
The final row describes the case when $x=0$, $y=\infty$ and $z=1$.
\[
  \begin{array}{c|c|c}
  U_z   & U_x       & U_y \\ \hline
  WLR^k & WLR^{k-1} & W \\
  WRL^k & W         & WRL^{k-1} \\
  R^k   & R^{k-1}   & \infty \\
  L^k   & 0         & L^{k-1} \\
  \varnothing & 0         & \infty
  \end{array}
\]
In the table, $W$ denotes an arbitrary (possibly empty) word 
in the letters $L$ and $R$, and $k\in\bbZ_{\ge1}$. 
Also, $0,\infty$ in the two right columns stand for 
the boundary fractions $0/1,1/0$ themselves.
We can find some examples of these four cases 
in \cref{fig:w:SBtree} of the Stern--Brocot tree.  

For a word $W$, let $\sfP_W$ be the product of the $2\times2$ matrices 
attached to its letters as in \eqref{eq:w:sfV}, 
with $\sfP_{\varnothing}=I$ the identity matrix.
Then the elliptic Stern--Brocot column is $\sfV(W)=\sfP_W\svc{1}{1}$.
The weight attached to a letter is determined by the preceding prefix,
so the matrices in a common prefix agree for $U_x$, $U_y$, and $U_z$.

First, we consider the case described in the first row of the table.
For $U_z=WLR^k$, let $\gamma_0$ denote the weight attached to 
the letter $L$ immediately following $W$, 
and let $\gamma_i$ denote the weight attached to 
the $i$th subsequent letter $R$ for $1\le i\le k$.
Set
\begin{align}
  \sfP \ceq 
  \sfP_W \sfB_L(\gamma_0) \sfB_R(\gamma_1) \dotsm \sfB_R(\gamma_{k-1}).
\end{align}
Since $U_x=WLR^{k-1}$ and $U_y=W$,
the corresponding Stern--Brocot columns are given by
\begin{align}
  \sfV(z)=\sfP\sfB_R(\gamma_k)\svc{1}{1},\quad
  \sfV(x)=\sfP\svc{1}{1},\quad
  \sfV(y)=\sfP_W\svc{1}{1}.
\end{align}
Using the identities
\begin{align}
  \sfB_R(\gamma_k)\svc{1}{1}-\svc{1}{1}=\gamma_k\svc{1}{0},\quad
  \sfB_R(\gamma)\svc{1}{0}=\gamma\svc{1}{0},\quad
  \sfB_L(\gamma)\svc{1}{0}=\gamma\svc{1}{1},
\end{align}
we obtain
\begin{align}
  \sfV(z)-\sfV(x)
  =\gamma_0\gamma_1\dotsm\gamma_k\sfV(y)
  =\frac{\Pi(z)}{\Pi(y)}\sfV(y),
\end{align}
where the last equality follows from the definition of the accumulated weight
$\Pi$ and the fact that $U_z$ is obtained from $U_y=W$ by appending $LR^k$.
Hence \eqref{eq:F:Farey} is proved in this case.

Next, we consider the case described in the second row of the table.
For $U_z=WRL^k$, define $\gamma_0,\gamma_1,\dotsc,\gamma_k$ analogously and set
\begin{align}
  \sfP \ceq
  \sfP_W \sfB_R(\gamma_0) \sfB_L(\gamma_1) \dotsm \sfB_L(\gamma_{k-1}).
\end{align}
Since $U_x=W$ and $U_y=WRL^{k-1}$, we have
\begin{align}
  \sfV(z)=\sfP\sfB_L(\gamma_k)\svc{1}{1},\quad
  \sfV(x)=\sfP_W\svc{1}{1},\quad
  \sfV(y)=\sfP\svc{1}{1}.
\end{align}
Using the identities
\begin{align}
  \sfB_L(\gamma_k)\svc{1}{1}=\svc{0}{1}+\gamma_k\svc{1}{1},\quad
  \sfB_L(\gamma)\svc{0}{1}=\svc{0}{1},\quad
  \sfB_R(\gamma)\svc{0}{1}=\svc{1}{1},
\end{align}
we obtain $\sfP\svc{0}{1}=\sfV(x)$ and hence
\begin{align}
  \sfV(z)=\sfV(x)+\gamma_k\sfV(y).
\end{align}
Since $U_z=U_yL$, we have $\gamma_k=\Pi(z)/\Pi(y)$.
This proves \eqref{eq:F:Farey} in this case.

For $U_z=R^k$, its Farey parents correspond to $R^{k-1}$ and $\infty$.
The identities
\begin{align}
  \sfB_R(\gamma)\svc{1}{1}=\svc{1}{1}+\gamma\svc{1}{0},\quad
  \sfB_R(\gamma)\svc{0}{0}=\svc{0}{0},\quad
  \sfB_R(\gamma)\svc{1}{0}=\gamma\svc{1}{0}
\end{align}
yield
\begin{align}
  \sfV(R^k)=\sfV(R^{k-1})+\Pi(R^k)\sfV(\infty).
\end{align}

For $U_z=L^k$, its Farey parents correspond to $0$ and $L^{k-1}$.
Let $\gamma_k$ denote the weight attached to the final letter $L$.
Using the identities
\begin{align}
  \sfB_L(\gamma_k)\svc{1}{1}=\svc{0}{1}+\gamma_k\svc{1}{1},\quad
  \sfB_L(\gamma)\svc{0}{1}=\svc{0}{1},
\end{align}
we obtain
\begin{align}
  \sfV(L^k)=\sfV(0)+\gamma_k\sfV(L^{k-1}),\quad
  \gamma_k=\frac{\Pi(L^k)}{\Pi(L^{k-1})}.
\end{align}

Finally, for the initial mediant $1=0\oplus\infty$, 
we have $\svc{1}{1}=\svc{0}{1}+\svc{1}{0}$.
Since $\Pi(1)=\Pi(\infty)=1$, 
this establishes \eqref{eq:F:Farey} for the initial edge,
and hence completes the proof.
\end{proof}

\begin{rmk}
In the $q$-degeneration, 
the recursion formula \eqref{eq:F:Farey} in \cref{thm:F:Farey}
reduces to the weighted mediant rule underlying \cite[Theorem~3]{MO}.
\end{rmk}

\begin{cor}\label[cor]{cor:F:edge}
For any Farey neighbors $0\le x<y\le\infty$, we have
\begin{align}\label{eq:F:edge}
  \det\bigl(\sfV(y),\sfV(x)\bigr)=\Pi(y).
\end{align}
\end{cor}

\begin{proof}
For the initial edge $(0,\infty)$, we have
$\det(\sfV(\infty),\sfV(0))=1=\Pi(\infty)$.
Suppose that \eqref{eq:F:edge} holds for a pair of Farey neighbors $x<y$, 
and let $z=x\oplus y$ be their mediant.
By \eqref{eq:F:Farey} and the bilinearity of the determinant, we obtain
\begin{equation}
  \det(\sfV(y),\sfV(z))=\det(\sfV(y),\sfV(x)),\quad
  \det(\sfV(z),\sfV(x))=\frac{\Pi(z)}{\Pi(y)}\det(\sfV(y),\sfV(x)).
\end{equation}
Thus, \eqref{eq:F:edge} also holds for the pairs $(x,z)$ and $(z,y)$.
This completes the proof 
by induction on the successive subdivisions of the initial edge.
\end{proof}


\begin{eg}\label[eg]{eg:farey-check}
We verify \cref{thm:F:Farey} and \cref{cor:F:edge}
for the rational numbers considered in \cref{eg:w:SBcol}.
Using the weights $w_{st}$ introduced in \cref{eg:w:weg}, we have
\begin{align*}
  \sfV(1)=\vc{1}{1},\quad
  \sfV(3)=\sfB_R(w_{11})\sfB_R(w_{12})\vc{1}{1}=\vc{1+w_{11}+w_{11}w_{12}}{1}.
\end{align*}
The accumulated weights of the relevant rational numbers are given by
\begin{gather*}
  \Pi(1)=1,\quad
  \Pi(2)=w_{11},\quad
  \Pi(3)=\Pi(3/2)=w_{11}w_{12},\\
  \Pi(5/2)=w_{11}w_{12}w_{13},\quad
  \Pi(5/3)=w_{11}w_{12}w_{22}.
\end{gather*}

Applying \eqref{eq:F:Farey} to the Farey mediants
$\frac{3}{2}=1\oplus2$,
$\frac{5}{2}=2\oplus3$, and
$\frac{5}{3}=\frac{3}{2}\oplus2$,
we obtain the following three identities:
\begin{gather*}
  \sfV\!\left(\frac{3}{2}\right) = \sfV(1)+w_{12}\sfV(2), \quad
  w_{12}=\frac{\Pi(3/2)}{\Pi(2)},\\
  \sfV\!\left(\frac{5}{2}\right) = \sfV(2)+w_{13}\sfV(3), \quad
  w_{13}=\frac{\Pi(5/2)}{\Pi(3)},\\
  \sfV\!\left(\frac{5}{3}\right) =
  \sfV\!\left(\frac{3}{2}\right)+w_{12}w_{22}\sfV(2), \quad
  w_{12}w_{22}=\frac{\Pi(5/3)}{\Pi(2)}.
\end{gather*}
These identities can also be verified directly using 
the explicit expressions for the Stern--Brocot columns in \cref{eg:w:SBcol}.

Similarly, direct computation gives 
the following four instances of \eqref{eq:F:edge}:
\begin{align*}
  \det\bigl(\sfV(2),\sfV(1)\bigr)    &=w_{11}=\Pi(2), &
  \det\bigl(\sfV(3),\sfV(2)\bigr)    &=w_{11}w_{12}=\Pi(3),\\
  \det\bigl(\sfV(5/2),\sfV(2)\bigr)  &=w_{11}w_{12}w_{13}=\Pi(5/2), &
  \det\bigl(\sfV(5/3),\sfV(3/2)\bigr)&=w_{11}w_{12}w_{22}=\Pi(5/3).
\end{align*}
\end{eg}


Next, we study the effect of the involution $x\mto1/x$ of 
the Stern--Brocot tree on our elliptic rational numbers.
This involution corresponds to interchanging 
the letters $R$ and $L$ in Stern--Brocot words.
Although Schlosser's elliptic weights \eqref{eq:w:sew} 
are not symmetric in their two arguments, 
they satisfy the reciprocal symmetry expressed by \eqref{eq:F:winv} below,
in which the two arguments are interchanged and 
the parameters are transformed simultaneously.

Set
\begin{align}
  (a^\vee,b^\vee,q^\vee)\ceq(b^{-1},a^{-1},q^{-1}), \quad
  J\ceq\begin{psmallmatrix}0&1\\1&0\end{psmallmatrix}.
\end{align}
Then, for all $s,t\in\bbZ_{\ge1}$, we have
\begin{align}\label{eq:F:winv}
  w_{a,b;q,p}(s,t) \, w_{a^\vee,b^\vee;q^\vee,p}(t,s)=1.
\end{align}

\begin{prp}[Elliptic inversion]\label[prp]{prp:F:inv}
For each positive rational number $x$, 
we write $\Pi_{a,b;q,p}(x)$ for the accumulated weight \eqref{eq:F:Pi}, 
making its dependence on the parameters explicit.
Then
\begin{align}\label{eq:F:cinv}
  \sfV_{a,b;q,p}(1/x) = \Pi_{a,b;q,p}(1/x)\,J\sfV_{a^\vee,b^\vee;q^\vee,p}(x),
\end{align}
and consequently
\begin{align}\label{eq:F:Einv}
  \clE_{a,b;q,p}(1/x) = \frac{1}{\clE_{b^{-1},a^{-1};q^{-1},p}(x)}.
\end{align}
\end{prp}

\begin{proof}
Set
$A\ceq aq^{s+2t}$,
$B\ceq bq^{2s+t-2}$, and
$C\ceq(a/b)q^{t-s-1}$.
The defining formula \eqref{eq:w:sew} gives
\begin{align}
\label{eq:F:iwp1}
  w_{a,b;q,p}(s,t)
  &=q\frac{\theta(A,B,C;p)}{\theta(q^{-2}A,q^2B,q^2C;p)},\\
\label{eq:F:iwp2}
  w_{a^\vee,b^\vee;q^\vee,p}(t,s)
  &=q^{-1}\frac{\theta((q^2B)^{-1},(q^{-2}A)^{-1},q^2C;p)}
               {\theta(B^{-1},A^{-1},C;p)}.
\end{align}
Applying the identity $\theta(u^{-1};p)=-u^{-1}\theta(u;p)$
to the theta factors with inverted arguments in \eqref{eq:F:iwp2},
we see that the right-hand side of \eqref{eq:F:iwp2}
is the reciprocal of that of \eqref{eq:F:iwp1}.
This proves \eqref{eq:F:winv}.

Write $U_x=\ve_1\ve_2\dotsm\ve_N$, and 
let $\ol{\ve}_j$ denote the letter obtained from $\ve_j$ 
by interchanging $R$ and $L$.
Then the Stern--Brocot word of $1/x$ is
$\ol{\ve}_1\ol{\ve}_2\dotsm\ol{\ve}_N$.
Let $\gamma_j$ be the weight attached to $\ve_j$ in $U_x$
with parameters $(a^\vee,b^\vee,q^\vee,p)$.
Since the numbers of occurrences of $R$ and $L$ in each prefix
are interchanged under the involution,
\eqref{eq:F:winv} implies that the weight attached to $\ol{\ve}_j$
in $U_{1/x}$ with parameters $(a,b,q,p)$ is $\gamma_j^{-1}$.
A direct calculation shows that, for $\ve\in\{L,R\}$,
\begin{align}
  \sfB_{\ol{\ve}}(\gamma^{-1}) = \gamma^{-1}J\,\sfB_{\ve}(\gamma)\,J.
\end{align}
Applying this identity to each factor in the matrix product
defining $\sfV_{a,b;q,p}(1/x)$
and using $J^2=I$ and $J\svc{1}{1}=\svc{1}{1}$, we obtain
\begin{align}
  \sfV_{a,b;q,p}(1/x) =
  \bigl(\tprd_{j=1}^{N}\gamma_j^{-1}\bigr)
  J\,\sfV_{a^\vee,b^\vee;q^\vee,p}(x).
\end{align}
By the definition \eqref{eq:F:Pi} of the accumulated weight,
the scalar factor on the right-hand side is $\Pi_{a,b;q,p}(1/x)$.
This proves \eqref{eq:F:cinv}.
Taking the ratio of the first entry to the second in both columns
of \eqref{eq:F:cinv} gives \eqref{eq:F:Einv}.
\end{proof}

In the $q$-degeneration $p\to0$, $a\to0$ and $b\to0$,
the weights $w_{a,b;q,p}(s,t)$ tend to $q$.
It follows from \eqref{eq:F:winv} that
$w_{a^\vee,b^\vee;q^\vee,p}(t,s) =w_{a,b;q,p}(s,t)^{-1}\to q^{-1}$.
Thus, \eqref{eq:F:Einv} reduces to the inversion formula 
for $q$-rational numbers in \cite{MO}:
\begin{align}
  [1/x]_q=\frac{1}{[x]_{q^{-1}}}.
\end{align}

\begin{rmk}\label[rmk]{rmk:F:cmp}
The parameter transformations in \cref{prp:w:trs} and \cref{prp:F:inv}
are compatible with each other.
With $p$ fixed, we regard them as maps on the triples $(a,b;q)$ given by
\begin{align}
  \sigma_L(a,b;q) \ceq (aq^2,bq;q),\quad
  \sigma_R(a,b;q) \ceq (aq,bq^2;q),\quad
  \iota(a,b;q)    \ceq (b^{-1},a^{-1};q^{-1}).
\end{align}
Direct substitution shows that
\begin{align}
  \iota^2 = \mathrm{id},\quad
  \iota\circ\sigma_L = \sigma_R\circ\iota,\quad
  \iota\circ\sigma_R = \sigma_L\circ\iota.
\end{align}
Note that, in the compositions $\sigma_R\circ\iota$ and $\sigma_L\circ\iota$,
the maps $\sigma_R$ and $\sigma_L$ are applied to the triple $\iota(a,b;q)$,
whose third component is $q^{-1}$.

On the Stern--Brocot tree, the involution $x\mto1/x$
conjugates the transformation $x\mto x+1$ to the transformation
$x\mto x/(x+1)$, as follows from the identity
$1/(1/x+1)=x/(x+1)$.
Accordingly, the first identity in \eqref{eq:w:trsQ}
can be deduced from the second identity in \eqref{eq:w:trsQ}
using the inversion formula \eqref{eq:F:Einv}
and the above relations among the parameter transformations.
We leave the details to the reader.

Finally, let us consider the group generated by 
$\sigma_L$, $\sigma_R$, and $\iota$.
Since $\sigma_L$ and $\sigma_R$ commute and
\begin{align}
  \sigma_L^m \sigma_R^n (a,b;q) = (aq^{2m+n},bq^{m+2n};q),
\end{align}
they generate a free abelian group of rank two.
The involution $\iota$ interchanges the two generators by conjugation.
Thus, the group under consideration is isomorphic to the semidirect product
\begin{align}
  \bbZ^2 \rtimes \bbZ/2\bbZ,
\end{align}
where the nontrivial element of $\bbZ/2\bbZ$ acts on $\bbZ^2$ 
by interchanging the two components.
This is precisely the extended affine Weyl group of $\mathrm{GL}_2$.
\end{rmk}


We close this subsection by studying the universal positivity of our construction.
We replace each weight $\omega_j(U)=\wtf{1+\ell_j}{1+r_j}$ in \eqref{eq:w:omega}
with an indeterminate $X_{1+\ell_j,1+r_j}$.
Let $\bsX$ denote the family of independent commuting indeterminates
$X_{s,t}$ indexed by $s,t\in\bbZ_{\ge1}$.
For each $x\in\bbQ_{\ge0}\cup\{\infty\}$, we write the resulting Stern--Brocot column as
\begin{align}
  \sfV_{\bsX}(x) \ceq \vc{R_x(\bsX)}{S_x(\bsX)},
\end{align}
with the conventions $\sfV_{\bsX}(0)\ceq\svc{0}{1}$ and
$\sfV_{\bsX}(\infty)\ceq\svc{1}{0}$.
By construction, we have
\begin{align}\label{eq:F:srng}
  R_x(\bsX), S_x(\bsX) \in \bbZ_{\ge0}[\bsX],
\end{align}
since, after replacing the weights with the indeterminates $X_{s,t}$,
all entries of the matrices $\sfB_L$ and $\sfB_R$ in \eqref{eq:w:mat}
belong to the semiring $\bbZ_{\ge0}[\bsX]$.

\begin{prp}\label[prp]{prp:F:pos}
For any $x,y\in\bbQ_{\ge0}\cup\{\infty\}$ satisfying $y<x$, we have
\begin{align}\label{eq:F:pos}
  R_x(\bsX)S_y(\bsX)-S_x(\bsX)R_y(\bsX)
  \in \bbZ_{\ge0}[\bsX] \setminus\{0\}.
\end{align}
\end{prp}

\begin{proof}
Fix $0<y<x<\infty$, and let $W$ be the longest common prefix of 
their Stern--Brocot words $U_y$ and $U_x$.
Let $\sfP$ denote the matrix product associated with $W$.
Since $\det\sfB_L(\gamma)=\det\sfB_R(\gamma)=\gamma$,
the determinant $\det\sfP$ is a monomial in the variables $X_{s,t}$.
Set $\xi\ceq X_{1+\ell(W),1+r(W)}$, 
the variable attached to the first letter following $W$.

Suppose first that neither $U_y$ nor $U_x$ is a prefix of the other.
Since $y<x$, the first letter following $W$ is $L$ in $U_y$ and $R$ in $U_x$.
We can therefore write
\begin{align}
  \sfV_{\bsX}(y)=\sfP\,\sfB_L(\xi)\svc{A_1}{A_2}, \quad
  \sfV_{\bsX}(x)=\sfP\,\sfB_R(\xi)\svc{B_1}{B_2},
\end{align}
where the columns $\svc{A_1}{A_2}$ and $\svc{B_1}{B_2}$
are obtained by applying the remaining matrix factors to $\svc{1}{1}$.
Since each of the matrices $\sfB_L$ and $\sfB_R$ preserves the property 
that both entries of a column are nonzero polynomials 
with nonnegative integer coefficients, 
all four polynomials $A_1,A_2,B_1,B_2$ have this property.
A direct calculation gives
\begin{align}
  \det(\sfV_{\bsX}(x),\sfV_{\bsX}(y)) =
  (\det\sfP)\bigl(\xi^2B_1A_1+\xi B_1A_2+B_2A_2\bigr),
\end{align}
which is a nonzero polynomial with nonnegative integer coefficients.
This proves \eqref{eq:F:pos} in the present case.

If $U_y=W$, the first letter following $W$ in $U_x$ is $R$.
With the same notation as above, we have
\begin{align}
  \sfV_{\bsX}(y)=\sfP\svc{1}{1}, \quad
  \sfV_{\bsX}(x)=\sfP\sfB_R(\xi)\svc{B_1}{B_2},
\end{align}
and the determinant is $(\det\sfP)\xi B_1$.
Similarly, if $U_x=W$, the first letter following $W$ in $U_y$ is $L$, 
and the determinant is $(\det\sfP)A_2$.
In both cases, the determinant is a nonzero polynomial 
with nonnegative integer coefficients.

Finally, the boundary cases $x=\infty$ or $y=0$ follow from the identities
\begin{align}
  \det(\sfV_{\bsX}(x),\sfV_{\bsX}(0))=R_x(\bsX), \quad
  \det(\sfV_{\bsX}(\infty),\sfV_{\bsX}(y))=S_y(\bsX)
\end{align}
and \eqref{eq:F:srng}, together with the fact that 
$R_x(\bsX)$ is nonzero for $x>0$ and $S_y(\bsX)$ is nonzero for $y<\infty$.
This completes the proof.
\end{proof}

In the $q$-degeneration, \cref{prp:F:pos} recovers \cite[Theorem~2]{MO}.

\section{Elliptic analogue of the quiver closure formula}
\label{s:Q}

The aim of this section is to express 
the two entries of the elliptic Stern--Brocot column 
as explicit sums over subsets of the vertex sets of type $A$ quivers.
In \cref{thm:Q:ellF}, we identify these entries with the closure polynomials 
of the quivers associated with the corresponding Stern--Brocot word.
This extends the combinatorial formula of \cite[Theorem~4]{MO} 
to the elliptic setting.
In \cref{ss:Q:QG}, we explain the interpretation of the closure polynomials 
in terms of quiver Grassmannians and 
establish their basic properties (\cref{prp:Q:pos}).

\subsection{Type \texorpdfstring{$A$}{A} quivers and their weights}
\label{ss:Q:Q}

Let $x \in \bbQ_{>1}$ with regular continued fraction expansion 
$[a_1,\dotsc,a_{2m}]$, and set $N\ceq\sum_{i=1}^{2m}a_i-1$.
Following \cite[\S3]{MO}, 
we consider the quiver $Q_x$ of type $A_N$ with vertices $1,\dotsc,N$.
From left to right, the arrows are arranged in consecutive blocks 
with alternating orientations, starting with a block of leftward arrows.
The lengths of these blocks are
$a_1-1,a_2,a_3,\dotsc,a_{2m-1},a_{2m}-1$ if $m\ge2$,
and $a_1-1,a_2-1$ if $m=1$.

The case $m\ge2$ is illustrated schematically below.
A block of length zero is understood to contain no arrows.

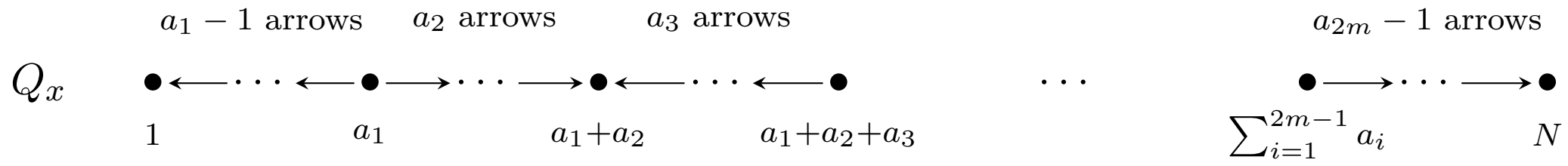
\begin{center}
\begin{tikzpicture}[>=stealth, every node/.style={inner sep=1pt}]
  \node      at (-1.00,0) {$Q_x$};
  \node (v0) at ( 0.00,0) {$\bullet$};
  \node (e1) at ( 0.95,0) {$\cdots$};
  \node (v1) at ( 1.90,0) {$\bullet$};
  \node (e2) at ( 2.90,0) {$\cdots$};
  \node (v2) at ( 3.90,0) {$\bullet$};
  \node (e3) at ( 4.95,0) {$\cdots$};
  \node (v3) at ( 6.00,0) {$\bullet$};
  \node      at ( 8.00,0) {$\cdots$};
  \node (v4) at (10.10,0) {$\bullet$};
  \node (e4) at (11.15,0) {$\cdots$};
  \node (v5) at (12.20,0) {$\bullet$};
  \draw[<-] (v0) -- (e1);
  \draw[<-] (e1) -- (v1);
  \draw[->] (v1) -- (e2);
  \draw[->] (e2) -- (v2);
  \draw[<-] (v2) -- (e3);
  \draw[<-] (e3) -- (v3);
  \draw[->] (v4) -- (e4);
  \draw[->] (e4) -- (v5);
  \node[font=\scriptsize] at ( 0.95,0.54) {$a_1-1$ arrows};
  \node[font=\scriptsize] at ( 2.90,0.54) {$a_2$ arrows};
  \node[font=\scriptsize] at ( 4.95,0.54) {$a_3$ arrows};
  \node[font=\scriptsize] at (11.15,0.54) {$a_{2m}-1$ arrows};
  \node at ( 0.0,-0.45) {$\scriptstyle 1$};
  \node at ( 1.9,-0.45) {$\scriptstyle a_1$};
  \node at ( 3.9,-0.45) {$\scriptstyle a_1+a_2$};
  \node at ( 6.0,-0.45) {$\scriptstyle a_1+a_2+a_3$};
  \node at (10.1,-0.45) {$\scriptstyle \sum_{i=1}^{2m-1}a_i$};
  \node at (12.2,-0.45) {$\scriptstyle N$};
\end{tikzpicture}
\end{center}

Write the Stern--Brocot word of $x$ as $U_x=\ve_1\ve_2\cdots\ve_N$.
Since $x>1$, we have $\ve_1=R$.
We associate vertex $i$ with the letter $\ve_i$ and assign to it the weight
\begin{align}\label{eq:Q:vwt}
  y_i \ceq \omega_i(U_x) = w_{a,b;q,p}(1+\ell_i,1+r_i) 
  \quad (i=1,2,\dotsc,N).
\end{align}
For $2\le i\le N$, the arrow between vertices $i-1$ and $i$ 
points to the left ($i-1\ot i$) if $\ve_i=R$ 
and to the right ($i-1\to i$) if $\ve_i=L$.
By \eqref{eq:cf:URL}, this rule agrees with the description of 
the blocks of arrows given above.

Following \cite[\S3.3]{MO}, we call a subset $C$ of the vertex set %
$\{1,2,\dotsc,N\}$ of a quiver $Q$ a \emph{closure} 
if no arrow of $Q$ has its source in $C$ and its target outside $C$.
Let $\Cl(Q)$ denote the set of all closures of $Q$.
For a quiver $Q$ with vertices $1,\dotsc,N$ and %
vertex weights $y_1,\dotsc,y_N$, we define its \emph{closure polynomial} by
\begin{align}\label{eq:Q:cls}
  F_Q(y_1,y_2,\dotsc,y_N) \ceq \sum_{C\in\Cl(Q)}\prod_{i\in C}y_i.
\end{align}
The closure polynomial of the empty quiver is $1$.

Let $Q'_x$ be the full subquiver of $Q_x$ on the vertex set 
$\{a_1+1,\dotsc,N\}$, with the vertex weights inherited from $Q_x$.
The quiver $Q'_x$ is empty when $x$ is an integer.
The case $m\ge2$ is illustrated schematically below.

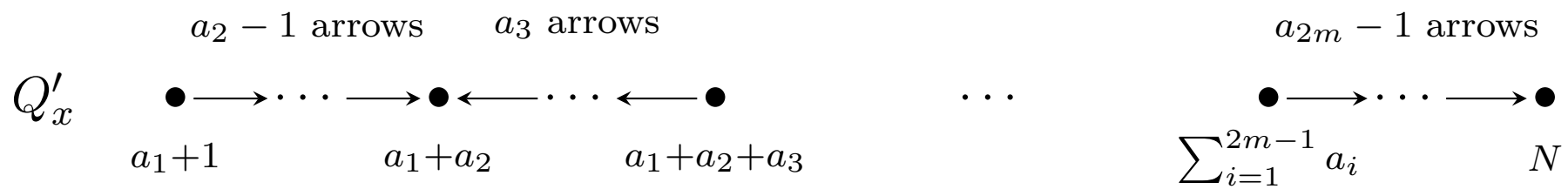
\begin{center}
\begin{tikzpicture}[>=stealth, every node/.style={inner sep=1pt}]
  \node      at ( 1.00,0) {$Q_x'$};
  \node (v0) at ( 2.00,0) {$\bullet$};
  \node (e1) at ( 3.00,0) {$\cdots$};
  \node (v1) at ( 4.00,0) {$\bullet$};
  \node (e2) at ( 5.05,0) {$\cdots$};
  \node (v2) at ( 6.10,0) {$\bullet$};
  \node      at ( 8.20,0) {$\cdots$};
  \node (v3) at (10.30,0) {$\bullet$};
  \node (e3) at (11.35,0) {$\cdots$};
  \node (v4) at (12.40,0) {$\bullet$};
  \draw[->] (v0) -- (e1);
  \draw[->] (e1) -- (v1);
  \draw[<-] (v1) -- (e2);
  \draw[<-] (e2) -- (v2);
  \draw[->] (v3) -- (e3);
  \draw[->] (e3) -- (v4);
  \node[font=\scriptsize] at ( 3.00,0.54) {$a_2-1$ arrows};
  \node[font=\scriptsize] at ( 5.05,0.54) {$a_3$ arrows};
  \node[font=\scriptsize] at (11.35,0.54) {$a_{2m}-1$ arrows};
  \node at ( 2.0,-0.45) {$\scriptstyle a_1+1$};
  \node at ( 4.0,-0.45) {$\scriptstyle a_1+a_2$};
  \node at ( 6.1,-0.45) {$\scriptstyle a_1+a_2+a_3$};
  \node at (10.3,-0.45) {$\scriptstyle \sum_{i=1}^{2m-1}a_i$};
  \node at (12.4,-0.45) {$\scriptstyle N$};
\end{tikzpicture}
\end{center}

\begin{thm}\label[thm]{thm:Q:ellF}
For every $x\in\bbQ_{>1}$,
\begin{align}\label{eq:Q:ellF}
  R_x=F_{Q_x}(y_1,\dotsc,y_N), \quad
  S_x=F_{Q'_x}(y_{a_1+1},\dotsc,y_N).
\end{align}
In the $q$-degeneration, this recovers \cite[Theorem~4]{MO}.
\end{thm}

\begin{proof}
For a subset $C$ of the vertex set $\{1,2,\dotsc,N\}$ of $Q_x$, 
define the binary variables $u_i$
by setting $u_i\ceq1$ if $i\in C$ and $u_i\ceq0$ otherwise.
For $i\ge2$, we call a pair $(u_{i-1},u_i)$ \emph{admissible} 
if it satisfies the closure condition 
for the arrow joining vertices $i-1$ and $i$.
If the arrow points to the left, $i-1\ot i$, the closure condition
$i\in C \Rightarrow i-1\in C$ is equivalent to $u_i\le u_{i-1}$.
Thus, the admissible pairs are $(1,1),(1,0),(0,0)$.
Similarly, if the arrow points to the right, $i-1\to i$, the closure condition
$i-1\in C\Rightarrow i\in C$ is equivalent to $u_{i-1}\le u_i$,
and the admissible pairs are $(1,1),(0,1),(0,0)$.

The \emph{transfer matrix} $\sfT_i$ encodes the admissibility condition
and incorporates the weight $y_i$ when vertex $i$ belongs to $C$.
Its rows and columns are indexed by the possible values of $u_{i-1}$ and $u_i$,
respectively, both in the order $(1,0)$.
The entries are defined by
\begin{align}
  (\sfT_i)_{\alpha,\beta}\ceq
  \begin{cases}
   y_i^\beta, &\text{$(\alpha,\beta)$ is admissible},\\
   0,         &\text{otherwise}.
  \end{cases}
\end{align}
It follows that, for $i\ge2$,
\begin{align}
  \sfT_i =
  \begin{cases}
   \begin{pmatrix}y_i&1\\   0&1\end{pmatrix}=\sfB_R(y_i), &i-1\ot i,\\[7pt]
   \begin{pmatrix}y_i&0\\ y_i&1\end{pmatrix}=\sfB_L(y_i), &i-1\to i.
  \end{cases}
\end{align}
For the first vertex, we introduce an auxiliary initial state
$u_0\in\{0,1\}$ and set $\sfT_1\ceq\sfB_R(y_1)$,
imposing the same admissibility condition $u_1\le u_0$ as for a leftward arrow.
Since $\ve_1=R$ and the orientation of the arrow between vertices $i-1$ and $i$
is determined by $\ve_i$ for $i\ge2$, we have
$\sfT_i=\sfB_{\ve_i}(y_i)$ for every $1\le i\le N$.
Hence, by \eqref{eq:w:sfV},
\begin{align}
  \sfT_1\dotsm\sfT_N\svc{1}{1}=\sfV(x).
\end{align}

Let us compute the two entries of the column $\sfT_1\dotsm\sfT_N\svc{1}{1}$ 
in terms of admissible sequences of states.
The entry corresponding to the initial state $u_0\in\{1,0\}$ is 
\begin{align}
  \sum_{u_1,\dotsc,u_N\in\{0,1\}}
  (\sfT_1)_{u_0,u_1} (\sfT_2)_{u_1,u_2} \dotsm (\sfT_N)_{u_{N-1},u_N}.
\end{align}
By the definition of $\sfT_i$, the summand vanishes unless every pair
$(u_{i-1},u_i)$ is admissible.
If all these pairs are admissible, the summand is equal to
$\prod_{i=1}^N y_i^{u_i}$.
Thus, the entry corresponding to $u_0$ is the sum of these monomials
over all admissible sequences $(u_1,\dotsc,u_N)$ 
with the prescribed initial state $u_0$.

Now we consider the entry corresponding to $u_0=1$.
In this case, both values of $u_1$ satisfy 
the initial admissibility condition $u_1\le u_0$.
There is therefore no additional restriction on the choice of the first vertex.
For $i\ge2$, the admissibility conditions are precisely 
the closure conditions for the arrows of $Q_x$.
Consequently, the admissible sequences correspond bijectively to 
the closures $C$ of $Q_x$, with $u_i=1$ if and only if $i\in C$.
Under this correspondence, 
the monomial $\prod_{i=1}^N y_i^{u_i}$ is equal to $\prod_{i\in C}y_i$.
Hence, the first entry is the closure polynomial of $Q_x$.

Next, consider the entry corresponding to $u_0=0$.
Suppose first that $x$ is not an integer.
The Stern--Brocot word $U_x$ begins with the block $R^{a_1}$, 
followed by the letter $L$.
The admissibility conditions for the initial block are
\begin{align}
  0=u_0\ge u_1\ge\dotsb\ge u_{a_1},
\end{align}
which force $u_1=\dotsb=u_{a_1}=0$.
The next letter is $L$, so the admissibility condition for the pair
$(u_{a_1},u_{a_1+1})$ is $u_{a_1}\le u_{a_1+1}$.
Since $u_{a_1}=0$, this condition imposes no restriction on $u_{a_1+1}$.
The remaining admissibility conditions are precisely 
the closure conditions for the arrows of the full subquiver $Q'_x$ 
on the vertex set $\{a_1+1,\dotsc,N\}$.
Thus, the admissible sequences with $u_0=0$ correspond 
bijectively to the closures of $Q'_x$.
Moreover, since the first $a_1$ states are zero, the corresponding monomial is
\begin{align}
  \prod_{i=1}^N y_i^{u_i} = \prod_{i=a_1+1}^N y_i^{u_i}.
\end{align}
Hence, the second entry is the closure polynomial of $Q'_x$.

If $x$ is an integer, all the letters in $U_x$ are $R$.
The initial condition $u_0=0$ then forces 
all the states $u_1,\dotsc,u_N$ to be zero.
In this case, the second entry is $1$, 
which agrees with the closure polynomial of the empty quiver $Q'_x$.
Now the proof is completed.
\end{proof}

\begin{eg}
For the rational number $5/2=[2,2]$, the Stern--Brocot word is $U_{5/2}=RRL$,
and the associated quiver is $Q_{5/2}=1\ot 2\to 3$,
with vertex weights $(y_1,y_2,y_3)=(w_{11},w_{12},w_{13})$.
The subquiver $Q'_{5/2}$ consists of the single vertex $3$ with no arrows.
The sets of closures of these quivers are
\begin{align}
  \Cl(Q_{ 5/2})=\{\varnothing,\{1\},\{3\},\{1,3\},\{1,2,3\}\},\quad
  \Cl(Q'_{5/2})=\{\varnothing,\{3\}\}.
\end{align}
Hence, by \cref{thm:Q:ellF}, we obtain
\begin{align}
  R_{5/2}=1+w_{11}+w_{13}+w_{11}w_{13}+w_{11}w_{12}w_{13}, \quad
  S_{5/2}=1+w_{13},
\end{align}
in agreement with the expressions in \cref{eg:w:SBcol}.

Similarly, for $5/3=[1,1,1,1]$, the Stern--Brocot word is $U_{5/3}=RLR$,
and the associated quiver is $Q_{5/3}=1\to 2\ot 3$,
with vertex weights $(y_1,y_2,y_3)=(w_{11},w_{12},w_{22})$.
The subquiver is $Q'_{5/3}=2\ot 3$,
and the sets of closures are
\begin{align}
  \Cl(Q_{ 5/3})=\{\varnothing,\{2\},\{1,2\},\{2,3\},\{1,2,3\}\},\quad
  \Cl(Q'_{5/3})=\{\varnothing,\{2\},\{2,3\}\}.
\end{align}
Hence, we obtain
\begin{align}
  R_{5/3}=1+w_{12}+w_{11}w_{12}+w_{12}w_{22}+w_{11}w_{12}w_{22}, \quad
  S_{5/3}=1+w_{12}+w_{12}w_{22}.
\end{align}
In both examples, specializing all the weights to $q$ recovers
the numerator and denominator polynomials given in \cite[Example~1.2]{MO}.
\end{eg}

\subsection{Quiver Grassmannians and coprimality}
\label{ss:Q:QG}

For general geometric properties of \emph{quiver Grassmannians} 
associated with representations of acyclic quivers, we refer to \cite{CR}.

Following the argument in \cite[\S3.3]{MO}, 
let $M_Q$ be the thin representation of a quiver $Q$ obtained by 
assigning a copy of $\bbC$ to each vertex and the identity map to each arrow.
Let $Q_0$ denote the vertex set of $Q$.
For a dimension vector $\bse=(e_i)_{i\in Q_0}\in\{0,1\}^{Q_0}$, 
we put $\supp(\bse)\ceq\{i\in Q_0\mid e_i=1\}$.
The corresponding quiver Grassmannian is given by
\begin{align}\label{eq:Q:qp}
  \Gr_{\bse}(M_Q) =
  \begin{cases}
  \{\pt\},     &\supp(\bse)\in\Cl(Q),\\
  \varnothing, &\text{otherwise}.
  \end{cases}
\end{align}

In the nonempty case, 
this quiver Grassmannian consists of a single reduced point.
Recall from \cite[(1.6)]{DWZ} that the \emph{$F$-polynomial} 
of a quiver representation $M$ of $Q$ is defined by
\begin{align}
  F_M(\bsy) \ceq \sum_{\bse\in\bbZ_{\ge0}^{Q_0}}
  \chi\bigl(\Gr_{\bse}(M)\bigr) \prod_{i\in Q_0}y_i^{e_i},
\end{align}
where $\chi$ denotes the topological Euler characteristic
and $\bsy=(y_i)_{i\in Q_0}$.
It follows from \eqref{eq:Q:qp} that the $F$-polynomial of $M_Q$
coincides with the closure polynomial $F_Q$ defined in \eqref{eq:Q:cls}.
This identification was used in \cite[\S3.3]{MO}.

In the formal setting of \cref{prp:F:pos},
the vertex weights $y_1,\dotsc,y_N$ of $Q_x$ are replaced by the variables
$X_{1+\ell_i,1+r_i}$ for $1\le i\le N$.
These variables are pairwise distinct since the lattice points 
corresponding to successive prefixes of a word are pairwise distinct.
We may therefore regard $y_1,\dotsc,y_N$ 
as independent commuting indeterminates.

\begin{prp}\label[prp]{prp:Q:pos}
Every monomial occurring in either $F_{Q_x}$ or $F_{Q'_x}$ 
is square-free and has coefficient $1$.
Both polynomials have constant term $1$, 
and they are coprime in $\bbZ[y_1,\dotsc,y_N]$.
\end{prp}

\begin{proof}
The assertions concerning the monomials and the constant terms 
follow immediately from the definition of the closure polynomial.

To prove coprimality, let $\sfM\ceq\sfT_1\dotsm\sfT_N$ be the product of
the transfer matrices introduced in the proof of \cref{thm:Q:ellF}.
Since $\det\sfT_i=y_i$ for every $i$, we have $\det\sfM=\prod_{i=1}^N y_i$.
By the proof of \cref{thm:Q:ellF}, we have
\begin{align}
  \sfM\vc{1}{1}=\vc{F_{Q_x}}{F_{Q'_x}}.
\end{align}
Writing $\sfM\svc{1}{0}=\svc{A}{B}$, we obtain
\begin{align}
  F_{Q_x}B-F_{Q'_x}A
  =(\det\sfM)\det\begin{pmatrix}1&1\\1&0\end{pmatrix}
  =-\prod_{i=1}^N y_i.
\end{align}
Any common irreducible factor of $F_{Q_x}$ and $F_{Q'_x}$ must therefore 
divide one of the variables $y_i$.
However, neither polynomial is divisible by any of these variables, 
since both have constant term $1$.
Thus, $F_{Q_x}$ and $F_{Q'_x}$ are coprime.
\end{proof}

\section{Vertex weights as elliptic Euler expressions}
\label{s:G}

As mentioned in \cref{ss:intro:geom}, 
the results obtained so far are algebraic in nature.
It is natural to ask whether they admit a geometric interpretation, 
possibly in terms of elliptic cohomology.
In this section, we present a preliminary attempt in this direction.

First, in \cref{ss:G:vwt}, we express the vertex weights
$y_i=\omega_i(U_x)$ in \eqref{eq:Q:vwt}
as \emph{meromorphic elliptic Euler expressions}
associated with balanced virtual representations of a torus.
We refer to \cref{rmk:G:Thom} 
for an explanation of the terminology ``Euler expressions''.
Next, in \cref{ss:G:root}, we explain how the dependence of these weights 
on Stern--Brocot prefixes can be described in terms of 
the root datum of $\SL_3$.
Finally, in \cref{ss:G:QGr}, we express the closure polynomials 
as finite traces over the quiver Grassmannians introduced in \cref{ss:Q:QG}.

\subsection{Elliptic weights and virtual representations}
\label{ss:G:vwt}

We refer to \cite[\S2]{AO} for a brief review of equivariant elliptic 
cohomology and to \cite{GKV,Konno} for related topics.

Fix $0<|p|<1$, and let $E_p\ceq\bbC^\times/p^{\bbZ}$ be
the elliptic curve with elliptic nome $p$.
Let $T\ceq(\bbC^\times)^3$ be the torus 
with multiplicative coordinate characters $Y,Z,q$.
The $T$-equivariant elliptic cohomology scheme of a point is given by
$\Ell_T(\pt)\simeq E_p^3$.
We consider the following virtual $T$-representation 
in the representation ring $R(T)=K_T(\pt)$:
\begin{align}\label{eq:G:vir}
  \clK(Y,Z;q) \ceq
  \bbC_Y+\bbC_Z+\bbC_{q^3YZ}-\bbC_{qYZ}-\bbC_{q^2Y}-\bbC_{q^2Z} = 
  (1-q^2)(Y+Z-qYZ).
\end{align}
Here we identify each character with 
the corresponding one-dimensional representation of $T$.

Using the modified theta function
\begin{align}\label{eq:G:vartheta}
  \vth{u} \ceq u^{-1/2}\tht{u},
\end{align}
we define the \emph{elliptic Euler expression} associated with $\clK(Y,Z;q)$ by
\begin{align}\label{eq:G:Euler}
  \eell(\clK(Y,Z;q)) \ceq
  \frac{\vth{Y}}{\vth{q^2Y}}
  \frac{\vth{Z}}{\vth{q^2Z}}
  \frac{\vth{q^3YZ}}{\vth{qYZ}} =
 q\frac{\tht{Y,Z,q^3YZ}}{\tht{q^2Y,q^2Z,qYZ}}.
\end{align}
The product of the characters occurring with positive signs in \eqref{eq:G:vir}
divided by the product of those occurring with negative signs is $q^{-2}$.
This accounts for the prefactor $q$ 
arising from the normalization factors $u^{-1/2}$.
Although the individual factors $\vth{u}$ require local choices of square roots,
the right-hand side of \eqref{eq:G:Euler} is 
a single-valued meromorphic function on $T$.
We therefore take this right-hand side 
as the definition of $\eell(\clK(Y,Z;q))$.

\begin{prp}\label[prp]{prp:G:Euler}
\begin{enumerate}
\item 
The elliptic Euler expression \eqref{eq:G:Euler} is invariant under
the transformations
$(Y,Z,q)\mto(p^{k_1}Y,p^{k_2}Z,p^{k_3}q)$
for all $k_1,k_2,k_3\in\bbZ$.
It therefore descends to a meromorphic function on $\Ell_T(\pt)\simeq E_p^3$.

\item
For $s,t\in\bbZ_{\ge1}$, set
\begin{align}\label{eq:G:YZ}
  Y \ceq bq^{2s+t-2},\quad Z \ceq (a/b)q^{t-s-1}.
\end{align}
Then Schlosser's small elliptic weight is given by
\begin{align}\label{eq:G:wYZ}
  w_{a,b;q,p}(s,t)=\eell(\clK(Y,Z;q)).
\end{align}
\end{enumerate}
\end{prp}

\begin{proof}
(1)
Let us rewrite the right-hand side of \eqref{eq:G:Euler} as
$q\prod_{j=1}^6\tht{u_j}^{\ep_j}$,
where $(u_1,u_2,u_3)=(Y,Z,q^3YZ)$ with $\ep_1=\ep_2=\ep_3=1$,
and $(u_4,u_5,u_6)=(q^2Y,q^2Z,qYZ)$ with $\ep_4=\ep_5=\ep_6=-1$.
Let $\bsm_j\in\bbZ^3$ denote the exponent vector of the monomial $u_j$
in the variables $Y,Z,q$, regarded as a column vector.
Explicitly, we have
\begin{gather*}
  \bsm_1=(1,0,0)^\top,\quad
  \bsm_2=(0,1,0)^\top,\quad
  \bsm_3=(1,1,3)^\top,\\
  \bsm_4=(1,0,2)^\top,\quad
  \bsm_5=(0,1,2)^\top,\quad
  \bsm_6=(1,1,1)^\top,
\end{gather*}
where $\top$ denotes the transpose.
A direct calculation shows that
\begin{align}\label{eq:G:anom}
  \sum_{j=1}^{6}\ep_j\bsm_j\bsm_j^{\top}=0,\quad
  \sum_{j=1}^{6}\ep_j\bsm_j=(0,0,-2)^\top.
\end{align}

Iterating the identity $\tht{pu}=-u^{-1}\tht{u}$, we obtain
\begin{align}
  \tht{p^nu}=(-1)^np^{-n(n-1)/2}u^{-n}\tht{u}
  \quad (n\in\bbZ).
\end{align}
Under the transformation
$(Y,Z,q)\mto(p^{k_1}Y,p^{k_2}Z,p^{k_3}q)$
associated with $\bsk\ceq(k_1,k_2,k_3)^\top$,
each monomial $u_j$ is multiplied by $p^{n_j}$,
where $n_j\ceq\bsm_j\cdot\bsk$.
Thus, the product of theta factors in \eqref{eq:G:Euler} is multiplied by
\begin{align}
  (-1)^{\sum_j\ep_jn_j}\,
  p^{-\frac12\sum_j\ep_jn_j^2+\frac12\sum_j\ep_jn_j}\,
  \prod_{j=1}^6u_j^{-\ep_jn_j}.
\end{align}
It follows from \eqref{eq:G:anom} that
\begin{align}
  \tsum_{j=1}^6\ep_jn_j^2 = 
  \bsk^{\top} \bigl(\tsum_{j=1}^6\ep_j\bsm_j\bsm_j^{\top}\bigr)\bsk = 0, 
  \quad
  \tsum_{j=1}^6\ep_jn_j = 
  \bigl(\tsum_{j=1}^6\ep_j\bsm_j\bigr)^\top\bsk = -2k_3.
\end{align}
Moreover, the exponent vector of the monomial
$\prod_{j=1}^6u_j^{-\ep_jn_j}$ in the variables $Y,Z,q$ is
\begin{align}
  -\tsum_{j=1}^6\ep_jn_j\bsm_j =
  -\bigl(\tsum_{j=1}^6\ep_j\bsm_j\bsm_j^{\top}\bigr)\bsk = 0.
\end{align}
Hence, this monomial is equal to $1$.
Combining these identities, we see that 
the product of theta factors is multiplied by $p^{-k_3}$.
This factor is exactly canceled by the factor $p^{k_3}$
arising from the transformation $q\mto p^{k_3}q$ of the prefactor.
Consequently, the elliptic Euler expression is invariant 
under all the prescribed transformations and 
descends to a meromorphic function on $E_p^3$.

(2)
Substituting the expressions for $Y$ and $Z$ in \eqref{eq:G:YZ}
into the right-hand side of \eqref{eq:G:Euler}
and using $YZ=aq^{s+2t-3}$,
we recover the defining expression \eqref{eq:w:sew}
for Schlosser's small elliptic weight.
This proves \eqref{eq:G:wYZ}.
\end{proof}

\begin{rmk}\label[rmk]{rmk:G:Thom}
In equivariant elliptic cohomology, 
the Euler class of a virtual representation is naturally regarded as 
a meromorphic section of the associated Thom line bundle over $\Ell_T(\pt)$.
For the virtual representation $\clK(Y,Z;q)$, 
\cref{prp:G:Euler} shows that the factors of automorphy cancel 
after the normalization in \eqref{eq:G:Euler}, 
yielding a meromorphic function on $E_p^3$.
We use the term elliptic Euler \emph{expression} 
to distinguish this explicitly normalized meromorphic function 
from the elliptic Euler \emph{class}, 
which is regarded as a section of a Thom line bundle.
The construction of the finite traces in \cref{ss:G:QGr} involves 
only these meromorphic functions.
It therefore requires neither a canonical trivialization of 
the Thom line bundles nor a pushforward in equivariant elliptic cohomology.
\end{rmk}

\subsection{Vertex characters and the root datum of %
 \texorpdfstring{$\SL_3$}{SL(3)}}
\label{ss:G:root}

By \eqref{eq:G:YZ}, the exponents of $q$ in the characters $Y$ and $Z$
associated with a lattice point $(s,t)$ are 
affine linear functions of $s$ and $t$.
Appending the letter $R$ corresponds to the shift $t \mto t+1$,
which multiplies $Y$ and $Z$ by $q$ and $q$, respectively.
Similarly, appending the letter $L$ corresponds to the shift $s \mto s+1$,
which multiplies $Y$ and $Z$ by $q^2$ and $q^{-1}$, respectively.
These shifts can be described in terms of the root datum of $\SL_3$.
This description provides a convenient way to express 
the changes in the exponents in terms of root-coweight pairings.

Let $H \subset \SL_3$ be the diagonal maximal torus,
$B \subset \SL_3$ the upper triangular Borel subgroup,
and $\frn$ the nilpotent radical of the Lie algebra of $B$.
Let $\alpha_1,\alpha_2$ be the simple roots determined by $B$, 
and set $Y \ceq e^{\alpha_1}$ and $Z\ceq e^{\alpha_2}$.
The $B$-representations
\begin{align}
  \frn/[\frn,\frn]\simeq\bbC_{\alpha_1}\oplus\bbC_{\alpha_2},\quad
  [\frn,\frn]\simeq\bbC_{\alpha_1+\alpha_2}
\end{align}
have $H$-characters $Y+Z$ and $YZ$, respectively.
Consequently,
\begin{align}
  (1-q^2)\bigl(\frn/[\frn,\frn]-q[\frn,\frn]\bigr)
\end{align}
defines a virtual representation of $B\times\bbC_q^\times$
whose character upon restriction to $H\times\bbC_q^\times$
is $\clK(Y,Z;q)$.

Define two cocharacters of $H$ by
\begin{align}
  \nu_R\ceq\alpha_1^\vee+\alpha_2^\vee=(1,0,-1),\quad
  \nu_L\ceq\alpha_1^\vee=(1,-1,0).
\end{align}
The pairings of the simple roots with these cocharacters are given by
\begin{align}\label{eq:G:pair}
\begin{array}{c|cc}
           & \nu_R & \nu_L \\ \hline
  \alpha_1 &    1  & 2 \\
  \alpha_2 &    1  & -1
\end{array}
\end{align}
Thus, multiplying each simple-root character $e^{\alpha_j}$ by
$q^{\langle\alpha_j,\nu_R\rangle}$ or
$q^{\langle\alpha_j,\nu_L\rangle}$ gives the transformations
\begin{align}\label{eq:G:Hecke}
  (Y,Z)\lmto(qY,qZ),\quad
  (Y,Z)\lmto(q^2Y,q^{-1}Z),
\end{align}
respectively.
These are precisely the transformations associated with 
appending the letters $R$ and $L$ to a Stern--Brocot prefix.

Now, let $x\in\bbQ_{>1}$, 
let $U_x=\ve_1\ve_2\dotsm\ve_N$ be its Stern--Brocot word, 
and let $r_i$ and $\ell_i$ be the numbers of occurrences of 
the letters $R$ and $L$, respectively, as defined in \cref{ss:w:SB}.
For each $i=1,2,\dotsc,N$, define the cocharacter
$\lambda_i\ceq r_i\nu_R+\ell_i\nu_L$.
With the initial values $Y_0\ceq bq$ and $Z_0\ceq a/(bq)$, we set
\begin{align}\label{eq:G:YL}
\begin{split}
 Y_i &\ceq Y_0q^{\langle\alpha_1,\lambda_i\rangle}
        = bq^{r_i+2\ell_i+1},\\
 Z_i &\ceq Z_0q^{\langle\alpha_2,\lambda_i\rangle}
        = (a/b)q^{r_i-\ell_i-1},\\
 \clK_i &\ceq \clK(Y_i,Z_i;q) \in R(T_{a,b,q}),
\end{split}
\end{align}
where $T_{a,b,q}\ceq(\bbC^\times)^3$ is 
the torus with coordinate characters $a,b,q$.

\begin{prp}\label[prp]{prp:G:rwt}
With the above notation, 
the vertex weights defined in \eqref{eq:Q:vwt} satisfy
\begin{align}\label{eq:G:mSch}
  y_i=\omega_i(U_x)=\eell(\clK_i).
\end{align}
Moreover, each vertex weight $y_i$ defines a meromorphic function on
$\Ell_{T_{a,b,q}}(\pt)\simeq E_p^3$.
\end{prp}

\begin{proof}
The explicit expressions for $Y_i$ and $Z_i$ in \eqref{eq:G:YL}
follow from the pairings in \eqref{eq:G:pair}.
Set $s\ceq1+\ell_i$ and $t\ceq1+r_i$.
Then we have $Y_i=bq^{2s+t-2}$ and $Z_i=(a/b)q^{t-s-1}$.
Thus, \eqref{eq:G:mSch} follows from \cref{prp:G:Euler}
and the definition of $\omega_i(U_x)$ in \eqref{eq:w:omega}.

Since $Y_i$ and $Z_i$ are Laurent monomials in $a,b,q$,
the assignment $(a,b,q)\mto(Y_i,Z_i,q)$ defines a homomorphism
of tori $T_{a,b,q}\to T$ for each $i$.
This homomorphism induces a morphism
\begin{align}
  \Ell_{T_{a,b,q}}(\pt)\to\Ell_T(\pt).
\end{align}
By \eqref{eq:G:mSch}, the vertex weight $y_i$ is the pullback of 
the meromorphic function $\eell(\clK(Y,Z;q))$ in \cref{prp:G:Euler} 
along this morphism.
This proves the last assertion.
\end{proof}

\begin{rmk}\label[rmk]{rmk:G:invK}
The elliptic inversion of \cref{prp:F:inv} can also be described 
in terms of the virtual representation \eqref{eq:G:vir}.
Consider the involutive automorphism of $T$ defined by
\begin{align}
  \iota_T(Y,Z,q)\ceq\bigl((qYZ)^{-1},q^2Z,q^{-1}\bigr).
\end{align}
Under the identification \eqref{eq:G:YZ}, this involution corresponds to
$(a,b;q;s,t)\mto(b^{-1},a^{-1};q^{-1};t,s)$.
Thus, \eqref{eq:F:winv} can be rewritten as
\begin{align}\label{eq:G:Einv}
  \eell\bigl(\clK(\iota_T(Y,Z,q))\bigr)\,
  \eell\bigl(\clK(Y,Z;q)\bigr)=1.
\end{align}

At the level of characters, we have
\begin{align}
  \iota_T^*\clK
  =q^{-1}Y^{-1}Z^{-1}+q^{-2}Y^{-1}+q^2Z-Y^{-1}-q^{-3}Y^{-1}Z^{-1}-Z.
\end{align}
Dualizing the four summands involving $Y^{-1}$ gives $-\clK(Y,Z;q)$.
Since $\vth{u^{-1}}=-\vth{u}$, these four dualizations leave the elliptic
Euler expression unchanged.
This provides another explanation of \eqref{eq:G:Einv}.

Moreover, $\iota_T$ interchanges the two character shifts in \eqref{eq:G:Hecke}
and sends the initial values $(bq,a/(bq))$ to
$(a^{-1}q^{-1},aq/b)$ for the transformed parameters.
These transformations are compatible with the exchange of the letters
$R$ and $L$ under $x\mto1/x$.
\end{rmk}

\subsection{Finite traces over quiver Grassmannians}
\label{ss:G:QGr}

Let $Q$ be a finite quiver of type $A$ with vertex set $Q_0$,
and let $M_Q$ be its thin representation with identity arrow maps,
as introduced in \cref{ss:Q:QG}.
We define the finite scheme
\begin{align}
  \frG_Q\ceq \coprod_{\bse\in\{0,1\}^{Q_0}}\Gr_{\bse}(M_Q).
\end{align}
By \eqref{eq:Q:qp}, this scheme is a disjoint union of reduced points
indexed by the closures $C\in\Cl(Q)$.
When $Q$ is empty, $\frG_Q$ consists of a single point.
For each $i\in Q_0$, let $\clU_i$ be the tautological vector bundle
on $\frG_Q$ associated with vertex $i$.
Its fiber at the point corresponding to $C\in\Cl(Q)$ is given by
\begin{align}\label{eq:G:taut}
  \clU_i\big|_C =
  \begin{cases}
    \bbC,&i\in C,\\
    0,&i\notin C.
  \end{cases}
\end{align}
We equip $\frG_Q$ with the trivial $T_{a,b,q}$-action
and each $\clU_i$ with the trivial equivariant structure.

Let $A\ceq\bbC(E_p^3)$ be the field of meromorphic functions on
$\Ell_{T_{a,b,q}}(\pt)\simeq E_p^3$.
We regard $\prod_{C\in\Cl(Q)}A$ as 
the ring of meromorphic coefficients on $\frG_Q$.
For the structure morphism $\pi_Q\colon\frG_Q\to\pt$,
we define the finite trace
\begin{align}
  (\pi_Q)_*\colon\prod_{C\in\Cl(Q)}A\to A,
  \quad (f_C)_C\lmto\sum_{C\in\Cl(Q)}f_C.
\end{align}

Now let $Q=Q_x$ be the quiver associated with $x\in\bbQ_{>1}$ 
introduced in \cref{ss:Q:Q}.
We attach to each vertex $i$ the virtual representation $\clK_i$
defined in \eqref{eq:G:YL}.
We define the elliptic Euler expression of the virtual bundle
$\bigoplus_i\clU_i\otimes\clK_i$ on $\frG_{Q_x}$
by specifying its value at each point $C\in\Cl(Q_x)$ as
\begin{align}\label{eq:G:UK}
  \rst{\eell\Bigl({\textstyle\bigoplus}_i\,\clU_i\otimes\clK_i\Bigr)}{C}
  \ceq\prod_{i\in C}\eell(\clK_i)\in A,
\end{align}
where the empty product is understood to be $1$.
This is the only extension of the notation $\eell$ to virtual bundles
on schemes that we use in this note.

\begin{prp}\label[prp]{prp:G:epush}
For $x\in\bbQ_{>1}$, let $\clK_i$ be the virtual representations 
defined in \eqref{eq:G:YL}, and consider the virtual bundle
\begin{align}\label{e:ellinsertion}
  \scK_x\ceq\bigoplus_{i=1}^{N}\clU_i\otimes\clK_i
  \quad\in K_{T_{a,b,q}}(\frG_{Q_x}).
\end{align}
Then we have
\begin{align}\label{e:ellpush}
  (\pi_{Q_x})_*\eell(\scK_x) = F_{Q_x}(y_1,\dotsc,y_N)=R_x.
\end{align}
For the full subquiver $Q'_x$ with the inherited vertex weights,
we define the virtual bundle $\scK'_x$ on $\frG_{Q'_x}$ analogously,
using only the indices $i=a_1+1,\dotsc,N$, where $a_1$ is 
the first entry of the regular continued fraction expansion of $x$.
Then we have
\begin{align}\label{e:ellpushprime}
  (\pi_{Q'_x})_*\eell(\scK'_x) = F_{Q'_x}(y_{a_1+1},\dotsc,y_N)=S_x.
\end{align}
\end{prp}

\begin{proof}
At the point of $\frG_{Q_x}$ corresponding to a closure $C\in\Cl(Q_x)$,
it follows from \eqref{eq:G:taut} that
$\scK_x|_C=\bigoplus_{i\in C}\clK_i$.
By \eqref{eq:G:UK} and \eqref{eq:G:mSch}, we therefore have
\begin{align}
  \rst{\eell(\scK_x)}{C}
  =\prod_{i\in C}\eell(\clK_i)
  =\prod_{i\in C}y_i.
\end{align}
By the definition of the finite trace, we obtain
\begin{align}
  (\pi_{Q_x})_*\eell(\scK_x)
  =\sum_{C\in\Cl(Q_x)}\prod_{i\in C}y_i
  =F_{Q_x}(y_1,\dotsc,y_N),
\end{align}
where the last equality follows from the definition of 
the closure polynomial in \eqref{eq:Q:cls}.
By \cref{thm:Q:ellF}, this closure polynomial is equal to $R_x$,
which proves \eqref{e:ellpush}.

The proof of \eqref{e:ellpushprime} is analogous.
When $Q'_x$ is empty, the scheme $\frG_{Q'_x}$ consists of a single point,
and all three expressions in \eqref{e:ellpushprime} are equal to $1$.
\end{proof}

\section{Further directions}
\label{s:conc}

\Cref{prp:G:epush} expresses the elliptic specializations of 
the closure polynomials in terms of equivariant data.
A natural question is whether the virtual representation $(1-q^2)(Y+Z-qYZ)$, 
including its $q$-grading, can be constructed canonically from 
the ungraded thin representation $M_{Q_x}$.
Indeed, since $Q_x$ is connected, we have
$\operatorname{Aut}(M_{Q_x})\simeq\bbC^\times$,
whereas the characters $Y_i,Z_i$ in \eqref{eq:G:YL} 
vary from vertex to vertex when $N\ge2$.
Thus, the character shifts along the vertices require additional data 
beyond the natural automorphism torus of $M_{Q_x}$.
Possible approaches to an intrinsic construction include the use of 
graded quiver representations or distinguished equivariant structures.
It would also be interesting to explore a connection with 
elliptic stable envelopes of quiver varieties \cite{AO,Konno}.

On the algebraic side, \cref{prp:w:trs} and \cref{prp:F:inv} describe 
the behavior of elliptic rational numbers 
under the three elementary transformations
$x\mto x+1$, $x\mto x/(x+1)$, and $x\mto1/x$,
together with the corresponding transformations of the parameters.
The compatibility of these parameter transformations is 
discussed in \cref{rmk:F:cmp}.
In particular, the inversion formula involves 
interchanging and inverting $a$ and $b$ and replacing $q$ by $q^{-1}$.
As explained in \cref{rmk:G:invK}, 
the involution $\iota_T$ of $T$ transforms $\clK$ into $-\clK$ 
after dualizing certain one-dimensional summands and 
interchanges the two character shifts.
It remains to understand whether this symmetry can be interpreted intrinsically
in terms of a duality of the representation-theoretic data 
underlying the virtual representations $\clK_i$.

Another direction is to replace Schlosser's small elliptic weights 
with nonzero weights constructed from shifted elliptic binomial coefficients 
or partition functions for lattice paths in rectangles.
The word identities in \cref{s:cf}, the Farey recursions, 
determinant identities, and formal positivity results in \cref{s:F}, 
as well as the quiver closure formulas in \cref{s:Q}, 
remain valid provided that the weights are assigned 
consistently to Stern--Brocot prefixes.
A natural problem is to find a systematic construction of such weights 
from Schlosser's small elliptic weights and 
to investigate their possible geometric interpretations.

\subsection*{Acknowledgements}

The author thanks Valentin Ovsienko for his talk on \cite{MO}
at the workshop ``International PhD Lab in Mathematics,''
held at Universit\'{e} de Reims Champagne-Ardenne in Reims in June 2018,
which inspired the present note.

The author was supported by JSPS KAKENHI Grant Number 26K06762
and the Asahipen Hikari Foundation.

\subsection*{Use of AI tools}

This note was prepared with the assistance of 
OpenAI's ChatGPT (GPT-6, GPT-6.2) and Anthropic's Claude (Opus 5.5). 
ChatGPT assisted in the theoretical development, 
the construction of illustrative examples, 
the preparation of Ti$k$Z code for figures, 
the identification of relevant references,
and the revision of the manuscript. 
Claude was used in the final stages of manuscript preparation 
to assist with proofreading and suggest improvements to the exposition.

The author independently verified all arguments, computations, 
and references in the final manuscript 
and takes full responsibility for any remaining errors.


\end{document}